\newif\ifpdf
\ifx\pdfoutput\undefined
\pdffalse 
\else
\pdfoutput=1 
\pdftrue \fi

\documentclass[reqno,twoside,11pt]{amsart}
\usepackage{amsmath}
\usepackage{amsfonts}
\usepackage{amssymb}
\usepackage{verbatim}

\IfFileExists{myowntimes.sty}{\usepackage{myowntimes}}
{\usepackage{times}\usepackage{mathrsfs}}

\DeclareFontFamily{OT1}{fraktura}{}
\DeclareFontShape{OT1}{fraktura}{m}{n} {<5> <6> <7> <8> <9> <10>
<11> <12> <13> <14.4> [1.1] eufm10}{}
\DeclareMathAlphabet{\fraktura}{OT1}{fraktura}{m}{n}

\newtheoremstyle{thm}{1.5ex}{1.5ex}{\itshape\rmfamily}{}
{\bfseries\rmfamily}{}{2ex}{}

\newtheoremstyle{rem}{1.3ex}{1.3ex}{\rmfamily}{}
{\itshape}
{} {1.5ex}{}

\newenvironment{proofsect}[1]
{\vskip0.1cm\noindent{\rmfamily\itshape #1.}}{\qed\vglue0.2cm}

\theoremstyle{thm}
\newtheorem{theorem}{Theorem}[section]
\newtheorem{lemma}[theorem]{Lemma}
\newtheorem{proposition}[theorem]{Proposition}

\newtheorem*{Main Theorem}{Main Theorem.}
\newtheorem{corollary}[theorem]{Corollary}
\newtheorem*{definition}{Definition.}

\theoremstyle{rem}
\newtheorem{remark}{{\itshape Remark}}[]

\numberwithin{equation}{section}

\renewcommand{\section}{\secdef\sct\sect}
\newcommand{\sct}[2][default]{\refstepcounter{section}
\addcontentsline{toc}{section}
{{\tocsection {}{\thesection}{\!\!\!\!#1\dotfill}}{}}
\vspace{0.7cm}
\centerline{ 
\scshape\arabic{section}.\ #1} \nopagebreak \vspace{0.2cm}}
\newcommand{\sect}[1]{
\vspace{0.4cm} \centerline{\large\scshape\rmfamily #1}
\vspace{0.2cm}}

\renewcommand{\subsection}{\secdef\subsct\sbsect}
\newcommand{\subsct}[2][default]{\refstepcounter{subsection}
\addcontentsline{toc}{subsection}
{{\tocsection{\!\!}{\hspace{1.2em}\thesubsection}{\!\!\!\!#1\dotfill}}{}}
\nopagebreak \vspace{0.45\baselineskip} {\flushleft\bf
\arabic{section}.\arabic{subsection}~\bf #1.~}
\\*[3mm]\noindent
\nopagebreak}
\newcommand{\sbsect}[1]{\vspace{0.1cm}\noindent
\textbf{#1.~}\vspace{0.1cm}}

\renewcommand{\subsubsection}{%
\secdef \subsubsect\sbsbsect}
\newcommand{\subsubsect}[2][default]{%
\refstepcounter{subsubsection} 
\addcontentsline{toc}{subsubsection}{{\tocsection{\!\!}
{\hspace{3.05em}\thesubsubsection}{\!\!\!\!#1\dotfill}}{}}
\nopagebreak
\vspace{0.15\baselineskip} \nopagebreak {\flushleft\rmfamily
\itshape\arabic{section}.\arabic{subsection}.\arabic{subsubsection}
\ \rmfamily #1\/.}\ }
\newcommand{\sbsbsect}[1]{\vspace{0.1cm}\noindent
\rmfamily \itshape
\arabic{section}.\arabic{subsection}.\arabic{subsubsection} \
\sffamily #1\/.\ }

\newcommand{\esssup}{\operatorname{esssup}}

\newcommand{\textd}{\text{\rm d}}

\newcommand{\Prob}{\text{\rm Prob}}

\newcommand{\BB}{\mathcal B}

\newcommand{\DD}{\mathcal D}
\newcommand{\EE}{\mathcal E}

\newcommand{\GG}{\mathcal G}
\newcommand{\HH}{\mathcal H}

\newcommand{\MM}{\mathcal M}
\newcommand{\NN}{\mathcal N}

\newcommand{\TT}{\mathcal T}
\newcommand{\UU}{\mathcal U}

\newcommand{\E}{\mathbb E}

\newcommand{\G}{\mathbb G}
\newcommand{\BbbH}{\mathbb H}

\newcommand{\N}{\mathbb N}

\newcommand{\BbbP}{\mathbb P}

\newcommand{\R}{\mathbb R}
\newcommand{\BbbS}{\mathbb S}
\newcommand{\T}{\mathbb T}

\newcommand{\Z}{\mathbb Z}

\newcommand{\twoeqref}[2]{(\ref{#1}-\ref{#2})}

\newcommand{\boldX}{{\boldsymbol X}}
\newcommand{\1}{{\text{\sf 1}}}
\newcommand{\zz}{{\fraktura z}}
\newcommand{\zzc}{\zz_{\text{\rm c}}}

\newcommand{\wBbbP}{\widehat\BbbP}
\newcommand{\wE}{\widehat\E}
\newcommand{\MMsharp}{\MM^\flat}

\newcommand{\eusmA}{\mathscr{A}}
\newcommand{\eusmB}{\mathscr{B}}
\newcommand{\eusmC}{\mathscr{C}}
\newcommand{\eusmG}{\mathscr{G}}

\begin{document}

\title[A model of organized criticality]
{\fontsize{16}{20}\selectfont Phase transition and critical behavior\\
in~a~model of organized criticality}
\author[M.~Biskup, Ph.~Blanchard, L.~Chayes, D.~Gandolfo and
T.~Kr\"uger] {M.~Biskup,${}^1\,$ Ph.~Blanchard,${}^2\,$
L.~Chayes,${}^1\,$ D.~Gandolfo,${}^{3,4}\,$ T.~Kr\"uger${}^2$}
\maketitle

\thispagestyle{empty}

\vspace{-3mm} \centerline{${}^1$\textit{Department of Mathematics,
UCLA, Los Angeles, California, USA}} \centerline{${}^2$\textit{Department
of Theoretical Physics, University of Bielefeld, Bielefeld,
Germany}} \centerline{${}^3$\textit{Phymath, Department of
Mathematics, University of Toulon, Toulon, France}}
\centerline{${}^4$\textit{CPT/CNRS, Luminy, Marseille, France}}
\vspace{4mm}

\begin{quote}
\footnotesize 
\textbf{Abstract:} 
We study a model of ``organized''
criticality, where a single avalanche propagates through an
\textit{a priori} static (i.e., organized) sandpile configuration.
The latter is chosen according to an i.i.d. distribution from a
Borel probability measure~$\rho$ on~$[0,1]$. The avalanche
dynamics is driven by a standard toppling rule, however, we
simplify the geometry by placing the problem on a directed, rooted
tree. As our main result, we characterize which~$\rho$ are
critical in the sense that they do not admit an infinite avalanche
but exhibit a power-law decay of avalanche sizes. Our analysis
reveals close connections to directed site-percolation, both in
the characterization of criticality and in the values of the
critical exponents.
\end{quote}

\setcounter{tocdepth}{3}

\footnotesize

\begin{list}{}
{\setlength{\topsep}{0in}\setlength{\leftmargin}{0.34in}\setlength{\rightmargin}{0.5in}}
\item[]
\tableofcontents
\end{list}
\vspace{-1.4cm}

\begin{list}{}
{\setlength{\topsep}{0in}\setlength{\leftmargin}{0.5in}\setlength{\rightmargin}{0.5in}}
\item[]
\hskip-0.01in
\hbox to 2.3cm{\hrulefill}
\item[]
{\fontsize{8.6}{8.6}\selectfont\copyright\,\,\,Copyright rests with the authors. Reproduction
of the entire article for non-commercial purposes
is permitted without charge.\vspace{2mm}}
\end{list}
\normalsize

\section{Introduction}

\subsection{Motivation}
\label{sec1.1}\noindent Since its discovery by Bak, Tang and
Wisenfeld~\cite{BTW1,BTW2}, self-organized criticality (SOC) has
received massive attention in the physics literature. Variants of
the original sandpile model of~\cite{BTW1} were studied and some
of them even ``exactly'' solved (see \cite{Dhar,Dhar-Ramaswamy,Maes1,Maes2} 
or \cite{Jensen} for a recent
review of the subject).  However, despite great efforts and
literally thousands of published papers, the present mathematical
understanding of SOC lags far behind the bold claims made by
physicists. Much of that failure can be attributed to the fact
that the models used to demonstrate SOC are difficult to formulate
precisely and/or too difficult to study using the current
techniques of probability theory and mathematical physics. From
the perspective of the latter fields, the situation seems ripe for
considering models which concern at least some aspects of SOC,
provided there is a decent prospect of a self-contained rigorous
analysis.

The general idea behind SOC models is very appealing. Consider for
instance Zhang's sandpile model~\cite{Zhang} on~$\Z^2$, where each
site has an energy variable which evolves in discrete time-steps
according to a simple ``toppling'' rule: If a variable exceeds a
threshold value, the excess is distributed equally among the
neighbors. The neighboring sites may thus turn supercritical and
the process continues until the excess is ``thrown overboard'' at
the system boundary. What makes this dynamical rule intriguing is
that, if the toppling is initiated from a ``highly excited''
state, then the terminal state (i.e., the state where the toppling
stops) is \textit{not} the most stable state, but one of many
\textit{least-stable}, stable states. Moreover, the latter state
is critical in the sense that further insertion of a small excess
typically leads to further large-scale events. Using the sandpile
analogy, such events are referred to as \textit{avalanches}.

In this paper, we study the scaling properties of a single
avalanche caused by an overflow at some site of a critical (i.e.,
least-stable) state. However, as indicated above, the full problem
is way too hard and we have to resort to simplifications. Our
simplifications are twofold: First, we treat the energy variables
of the critical state as independent and, second, we consider the
model on a directed, rooted tree rather than~$\Z^2$. The first
assumption is fairly reasonable, at least on a coarse-grained
scale, because numerical results~\cite{Grassberger-Manna} suggest
a rather fast decay of spatial correlations in the critical
states. The second assumption will allow us to treat the
correlations between different branches of the avalanche as
conditionally independent, which will greatly facilitate the
analysis. Finally, the reduced geometry allows for the existence
of a natural monotonicity not apparent in the full-fledged model.

While placing the model on a tree simplifies the underlying
geometry, some complexity is retained due to the generality of the
single-site energy variable distribution. In fact, the set of
underlying distributions plays the role of a parameter space in
our case. As our main result, we characterize the subspace of
distributions for which the configurations of energy variables
have exactly the behavior expected from the SOC states: \textit{no
infinite avalanches} but a \textit{power-law decay} of avalanche sizes. As
it turns out, there is a close connection to site-percolation on
the underlying graph, both in the characterization of criticality
and in the values of the critical exponents. However, the
significance of this connection for the general SOC models has not
yet been clarified.

\subsection{The model}
\label{sec1.2}\noindent In order to precisely define our
single-avalanche model, we need to introduce some notation. Let~$b>1$ be an integer and let~$\T_b$ be a~$b$-nary rooted tree, with
the root vertex denoted by~$\varnothing$. 
We use~$|\sigma|=k$ to denote that~$\sigma\in\T_b$ is on the~$k$-th
layer. When~$|\sigma|=k$, we represent~$\sigma$ as a~$k$-component object. Each
component is an integer in~$\{1,\dots,b\}$; hence the site label
can be used to trace the path from~$\sigma$ back to the root. If~$\sigma$ is an~$\ell$-th level site with~$\ell>0$, we let~$m(\sigma)$ denote the ``mother-site.'' Explicitly, if~$\sigma=(\sigma_1,\dots,\sigma_\ell)$, then~$m(\sigma)=(\sigma_1,\dots,\sigma_{\ell-1})$. The edges of~$\T_b$
are the usual directed edges~$\bigl\{(\sigma',\sigma)\in\T_b\times\T_b\colon
\sigma'=m(\sigma)\bigr\}$.

Let~$\MM$ be the space of all probability measures on the Borel~$\sigma$-algebra of~$[0,1]$. Fix a~$\rho\in\MM$ and let~$\BbbP_\rho=\rho^{\T_b}$. Let~$\E_\rho$ denote the expectation with
respect to~$\BbbP_\rho$. The dynamical rule driving the evolution
is  defined as follows: Let~$\boldX=(X_\sigma)_{\sigma\in\T_b}$ be
the collection of i.i.d.\ random variables with joint probability
distribution~$\BbbP_\rho$ and let~$v\in(0,\infty)$. The process
generates the sequence
\begin{equation}
\boldX^{(v)}(t)=\bigl(X_\sigma^{(v)}(t)\bigr)_{\sigma\in\T_b},
\qquad t=0,1,\dots,
\end{equation}
obtained from the initial condition
\begin{equation}
\label{initval} X_\sigma^{(v)}(0)=\begin{cases}
X_\varnothing+v, \qquad&\text{if }\sigma=\varnothing,\\
X_\sigma,\qquad&\text{otherwise,}
\end{cases}
\end{equation}
by successive applications of the deterministic (Markov) update
rule
\begin{equation}
\label{update} X_\sigma^{(v)}(t+1)=\begin{cases}
X_\sigma^{(v)}(t)+\frac1b X_{m(\sigma)}^{(v)}(t), \qquad&\text{if
} X_{m(\sigma)}(t)\ge1,
\\
0,\qquad&\text{if }X_\sigma^{(v)}(t)\ge1,\\
X_\sigma^{(v)}(t),\qquad&\text{otherwise.}
\end{cases}
\end{equation}
Note that, if~$X_\sigma^{(v)}(t+1)=X_\sigma^{(v)}(t)$ for all~$\sigma\in\T_b$, then~$X_\sigma^{(v)}(t)\le1$ and the process has
effectively \textit{stopped}. (However, we let~$X_\sigma^{(v)}(t)$
be defined by \eqref{update} for all~$t\ge0$.)

Here is an informal description of the above process: Starting at
the root we first check whether~$X_\varnothing+v\ge1$ or not. If
not, the process stops but if so, then this value is distributed
evenly among the ``daughter'' cells, which have their values
updated to~$X_\sigma^{(v)}(1)=X_\sigma+\frac1b(X_\varnothing+v)$.
The value~$X_\varnothing^{(v)}(1)$ is set to zero and we say that
the root has ``avalanched.'' If none of the updated ``daughter''
values exceed one, the process terminates; however, if there is
any first-level~$\sigma$ with~$X_\sigma^{(v)}(1)\ge1$, then~$X_\sigma^{(v)}(1)$ is set to zero, the value~$X_\sigma^{(v)}(1)$
is evenly distributed among the ``daughters'' of~$\sigma$ and we
say that~$\sigma$ has ``avalanched.'' The process at future times
is described similarly.

Obviously, the variables~$X_\sigma$ play the role of the ``energy
variables'' in the description of Zhang's avalanche model in
Section~\ref{sec1.1}. In our case the critical threshold is one,
but, in~\eqref{update}, we chose to distribute the entire value of
an ``avalanching'' site rather than just the excess to the
(forward) neighbors. This choice is slightly more advantageous
technically.

\subsection{Main questions and outline}
\label{sec1.3}\noindent Let~$\eusmA^{(v)}(t)=\{\sigma\in\T_b\colon
X_\sigma^{(v)}(t)=0,\,X_\sigma^{(v)}(s)\ne0\text{ for some }s<t\}$
be the set of sites that have ``avalanched'' by time~$t$.
Similarly, let~$\eusmA^{(v)}=\bigcup_{t\ge0}\eusmA^{(v)}(t)$ be
the set of sites that will ever avalanche. We use~$|\eusmA^{(v)}|$ to denote the number of sites in the avalanched set
(which includes the possibility of~$|\eusmA^{(v)}|=\infty$). The
set~$\eusmA^{(v)}$ and its dependence on~$\rho$ and~$v$ are the
primary focus of our study.

The first question is whether the process~$\boldX^{(v)}(t)$ lives
forever, i.e., is there an infinite avalanche? More precisely, for
what measures~$\rho\in\MM$ is the probability
\begin{equation} \label{Ainfty}
A_\infty^{(v)}=\BbbP_\rho\bigl(|\eusmA^{(v)}|=\infty\bigr)
\end{equation}
non-zero for some value of~$v$? A related question is whether the
average size of the avalanched set is finite. The relevant object
is defined~by
\begin{equation}
\label{chiv} \chi^{(v)}=\E_\rho\bigl(|\eusmA^{(v)}|\bigr).
\end{equation}
(Notice that, due to the directed nature of the dynamical rule,
both quantities~$A_\infty^{(v)}$ and~$\chi^{(v)}$ are monotone in
the underlying measure and~$v$.) Again, we ask: For what measures~$\rho$ we have~$\chi^{(v)}=\infty$ for some~$v$? In addition, we might
ask: Is the divergence of the mean avalanche size equivalent 
to the onset of infinite avalanches or can there be 
an~\textit{intermediate phase}?

To give answers to the above questions, we will
parametrize the set~$\MM$ by values of a particular functional~$\zz\colon\MM\to[0,1]$. Here~$\zz(\rho)$ roughly corresponds to
the conditional probability in distribution~$\BbbP_\rho$ that,
given the avalanche has reached a site~$\sigma\in\T_b$ far away
from the root, the site~$\sigma$ will also avalanche. (The
definition of~$\zz$ is somewhat technical and we refer the reader
to Section~\ref{sec2.2} for more details.) The characterization of
the avalanche regime in terms of~$\zz$ is then very transparent:
There is a \textit{critical value}~$\zzc=\frac1b$, such that the quantity~$\chi^{(v)}$ for measure~$\rho$ diverges if~$\zz(\rho)>\zzc$ and~$v$ is sufficiently large, while it is finite for all~$v$ if~$\zz(\rho)<\zzc$. Similarly we show, for a reduced class
of measures, that~$A_\infty^{(v)}$ for measure~$\rho$ vanishes for
all~$v$ if and only if~$\zz(\rho)\le\zzc$. These results are
formulated as Theorems~\ref{thm2.4} and~\ref{thm3.1} in
Sections~\ref{sec2.2} and~\ref{sec3.1}, respectively. (Outside the
reduced class of measures, there are some exceptions to the rule
that~$A_\infty^{(v)}\equiv0$ for measures~$\rho$ with~$\zz(\rho)=\zzc$,
i.e., there are some measures which avalanche also \textit{at} criticality,
see Remarks~\ref{rem1} and~\ref{rem2} in Section~\ref{sec2} for
more details. These examples are fairly contrived, so we exclude
them from further considerations.)

Note that, for both quantities \eqref{Ainfty} and \eqref{chiv},
the transitions happen at the same value,~$\zzc$, which rules out
the possibility of an intermediate phase. To elucidate the
behavior of~$\zz$ near~$\zzc$, it is worthwhile to introduce
appropriate \textit{critical exponents}. In particular, we ask whether
there is a critical exponent~$\gamma>0$ such that
\begin{equation}
\label{critexp1} \chi^{(v)}\sim
\bigl(\zzc-\zz(\rho)\bigr)^{-\gamma},\qquad \zz(\rho)\uparrow\zzc,
\end{equation}
an exponent~$\beta>0$ such that
\begin{equation}
\label{critexp2} A_\infty^{(v)}\sim
\bigl(\zz(\rho)-\zzc\bigr)^\beta,\qquad \zz(\rho)\downarrow\zzc,
\end{equation}
and, finally, an exponent~$\delta>0$ such that if~$\zz(\rho)=\zzc$, then
\begin{equation}
\label{critexp3} \BbbP_\rho\bigl(|\eusmA^{(v)}|\ge n\bigr)\sim
n^{-1/\delta},\qquad n\to\infty.
\end{equation}
All of these three relations of course include an appropriate
interpretation of the symbol ``$\sim$'' and, with the exception of
the last relation, also an interpretation of the limit
``$\zz(\rho)$ tends to~$\zzc$.''

The relations for the critical exponents are the subject of
Theorem~\ref{thm4.1} in Section~\ref{sec4}. The upshot is that all
three exponents take the \textit{mean-field percolation} values,
\begin{equation}
\gamma=1,\quad\beta=1,\quad\delta=2.
\end{equation}
Neither the fact that the critical value~$\zzc$ equals the
percolation threshold for site percolation on~$\T_b$ is a
coincidence. Indeed, the avalanche problem can be characterized in
terms of a correlated-percolation problem on~$\T_b$ (see
Section~\ref{sec2}).

We finish with a brief outline of the paper: Section~\ref{sec2}
contains our percolation criteria for the existence of infinite
avalanches leading naturally  to the definition of the functional~$\zz$. In Section~\ref{sec3} we show that~$\zzc=\frac1b$ is the
unique critical ``point'' of our model, thus ruling out the
possibility of an intermediate phase. Section~\ref{sec4} proves
the above relations for the critical exponents. Finally, in
Section~\ref{sec5} we develop a coupling argument which is the
core of the proofs  of the aforementioned results in
Sections~\ref{sec3} and~\ref{sec4}. The principal results of this paper are
Theorem~\ref{thm2.4} (Section~\ref{sec2.2}), Theorem~\ref{thm3.1}
(Section~\ref{sec3.1}) and Theorem~\ref{thm4.1}
(Section~\ref{sec4.1}).

\section{Percolation criteria}
\label{sec2}

\subsection{Simple percolation bounds}
\label{sec2.1}\noindent We start by deriving criteria for the
presence and absence of an infinite avalanche based on a
comparison to site percolation on~$\T_b$. Let~$x_\star$ denote the
maximum of the support of~$\rho$, i.e.,
\begin{equation}
\label{xstar} x_\star=\sup\bigl\{y\in[0,1]\colon
\rho([y,1])>0\bigr\},
\end{equation}
and let us define~$\theta_b$ by
\begin{equation}
\label{thetab} \theta_b=\frac b{b-1}x_\star.
\end{equation}
It is noted that if~$X_\varnothing+v\le\theta_b$, then the largest
value that~$X_\sigma(t)$ for any~$\sigma\in\T_b$ could conceivably
achieve (just prior to its own avalanche) is~$\theta_b$.

\smallskip
The following is based on straightforward percolation arguments:

\begin{proposition}
\label{prop2.1}
\settowidth{\leftmargini}{(11)}
\begin{enumerate}
\item[(1)]
If~$\rho([\frac{b-1}b,1])>\frac1b$, then~$\BbbP_\rho(|\eusmA^{(v)}|=\infty)>0$ for
all~$v>1-x_\star$.
\item[(2)]
If either~$\theta_b<1$ or~$\theta_b>1$ and~$\rho([1-\frac1b\theta_b,1])\le\frac1b$, then~$\BbbP_\rho(|\eusmA^{(v)}|=\infty)=0$
for all~$v\ge0$.
\end{enumerate}
\end{proposition}

In both cases we note that the quantity~$\frac1b$ on the
right-hand side of the inequalities is the percolation threshold
for~$\T_b$. Obviously, this is no coincidence; indeed, the proof of
part~(1) is easily generalizable to any transitive infinite graph.

\begin{proofsect}{Proof of Proposition~\ref{prop2.1}}
Let us start with (1): A site~$\sigma\ne\varnothing$ is called
occupied if~$X_\sigma\ge1-\frac1b$, while the root~$\varnothing$
is called occupied if~$X_\varnothing+v\ge1$. Denoting by~$\eusmC^{(v)}$ the connected component of occupied sites containing the origin, it is
not hard to see that~$\eusmA^{(v)}\supset \eusmC^{(v)}$. Indeed,
assuming~$X_\varnothing+v>1$, each daughter site of the origin
receives at least~$\frac1b$; those daughter sites~$\sigma$ with~$X_\sigma\ge 1-\frac1b$ will be triggered, which will in turn
cause avalanches in the next generation of occupied sites, etc.
Evidently, whenever the occupied sites percolate, there is an
infinite avalanche.

Part (2) is proved in a similar fashion. Suppose first that~$\theta_b>1$ and call a site~$\sigma\ne\varnothing$ occupied if~$X_\sigma\ge1-\frac1b\theta_b$, and vacant otherwise. The
definition is as before for the root. As observed previously, if~$X_\varnothing+v\le\theta_b$, then no site receives more than~$\frac1b\theta_b$ from its parent. Under these circumstances, a
vacant site will never avalanche and, denoting the occupied
cluster of the origin by~$\bar{\eusmC}^{(v)}$, we have~$\eusmA^{(v)}\subset\bar{\eusmC}^{(v)}$. Since~$\rho([1-\frac1b\theta_b,1])\le\frac1b$ was assumed, we have that~$|\bar{\eusmC}^{(v)}|<\infty$ almost surely and thus~$|\eusmA^{(v)}|<\infty$ whenever~$X_\varnothing+v\le\theta_b$. It is then
easy to show, however, that~$|\eusmA^{(v)}|<\infty$ almost surely
for all~$v\ge0$. Indeed, let~$k\ge0$ be an integer so large that
\begin{equation}
\label{kvrel} (x_\star+v-\theta_b)b^{-k}<\theta_b-1.
\end{equation}
If~$\sigma$ is a site with~$|\sigma|=k$ that has been reached by
an avalanche, then~$\sigma$ could not receive more~than
\begin{equation}
\label{asymptgift} x_\star\bigl(b^{-1}+\dots+b^{-(k-1)}\bigr)+
b^{-k}(X_\varnothing+v)=b^{-1}\theta_b+
b^{-k}(X_\varnothing+v-\theta_b)
\end{equation}
from its parent. Now, if~$\sigma$ is vacant, then the maximal
possible value for~$X_\sigma(k)$ (i.e., prior to its own
avalanche) is no larger than~$1+b^{-k}(x_\star+v-\theta_b)$. By
\eqref{kvrel}, this amount is strictly less than~$\theta_b$, so by
our previous reasoning,~$\sigma$ cannot give rise to an infinite
avalanche. By absence of percolation, there is a ``barrier''~$\BbbS_k$ of vacant sites above the~$(k+1)$-st layer in~$\T_b$,
that every path from the root to infinity must pass through. Our
previous arguments show that the avalanche cannot go beyond the
union of occupied connected components rooted at~$\BbbS_k$. Hence,~$|\eusmA^{(v)}|<\infty$ with probability one.

The case~$\theta_b<1$ is handled analogously. Indeed, a simple
calculation reveals that the right-hand side of \eqref{asymptgift}
plus~$x_\star$ is eventually strictly less than one and the
avalanche terminates within a deterministic ($v$-dependent) amount
of time.
\end{proofsect}

The arguments in the proof immediately give us the following
corollary:

\begin{corollary}
\label{cor2.2} If~$\theta_b\ne1$ and~$\rho([1-\frac1b\theta_b,\frac{b-1}b))=0$, then there is an
infinite avalanche if and only if occupied sites, i.e., sites~$\sigma\in\T_b$ with value~$X_\sigma\ge1-\frac1b$, percolate. In
addition, if~$X_\varnothing+v\le\theta_b$, then~$\eusmA^{(v)}$
coincides exactly with the occupied connected component
of~the~root.
\end{corollary}

\begin{remark}
\label{rem1} The exceptional cases,~$\theta_b=1$, can only arise
from the circumstance that~$x_\star=1-\frac1b$. (Notice that the
proof of Proposition~\ref{prop2.1}(2) does not apply because the
inequality \eqref{kvrel} cannot be satisfied.) For~$\theta_b=1$,
the situation is marginal and, in fact, slightly subtle. Indeed,
if~$x_\star=1-\frac1b$ and~$\BbbP(X\ge x_\star)=\frac1b$, then the
existence of an infinite avalanche depends on the detailed
asymptotic of~$\BbbP(X\ge x_\star-\epsilon)$ as~$\epsilon\downarrow0$, see Remark~\ref{rem2} in the next section.
We exclude the cases~$\theta_b=1$ from our analysis because we
believe that this ``pathological'' behavior is in no way~generic.
\end{remark}

\subsection{Phase transition}
\label{sec2.2}\noindent As is seen from Corollary~\ref{cor2.2}, in
certain cases the avalanche problem reduces to the usual
(independent) percolation model. The general problem can also be
presented as a percolation phenomenon albeit with correlations.
Indeed, let~$X_1,\dots,X_n$ are i.i.d.\ with distribution~$\rho$
and let
\begin{equation}
\label{Qntheta}
Q_n^{(\theta)}=X_n+\frac{X_{n-1}}b+\dots+\frac{X_1}{b^{n-1}}
+\frac\theta{b^n},
\end{equation}
In the case~$n=0$, we let~$Q_0^{(\theta)}=\theta$. Similarly, for~$\sigma\in\T_b$, we define~$Q_\sigma^{(\theta)}$ by
\eqref{Qntheta} with~$n=|\sigma|$ and~$X_1,\dots,X_{|\sigma|}$ being the values
along the unique path connecting~$\sigma$ to the root. Explicitly,
we set~$Q_\varnothing^{(\theta)}=\theta$ and define
\begin{equation}
\label{Qsigma} Q_\sigma^{(\theta)}=X_\sigma+\frac1b
Q_{m(\sigma)}^{(\theta)}, \qquad \sigma\ne\varnothing.
\end{equation}
Note that here~$\theta$ plays the role of the quantity~$X_\varnothing+v$. Clearly,~$Q_n^{(\theta)}\overset{\DD}=Q^{(\theta)}_\sigma$, whenever~$n=|\sigma|$.

\begin{proposition}
\label{prop2.3} Let~$v\ge0$ and let~$\theta=X_\varnothing+v$. For
each~$\sigma\in\T_b$, let us call~$\sigma$ \rm open \it if~$Q_\sigma^{(\theta)}\ge1$ and \rm closed \it otherwise. Then~$\sigma\in\eusmA^{(v)}$ if and only if~$\sigma$ belongs to the
open cluster containing the root. In particular, percolation of
open sites is the necessary and sufficient condition for infinite
avalanches.
\end{proposition}

\begin{proofsect}{Proof}
By definition,~$Q_\varnothing^{(\theta)}=\theta=X_\varnothing+v$. Now, if~$X_\sigma(t)=Q_\sigma^{(\theta)}$ for a site~$\sigma\in\T_b$ that
avalanches at time~$t=|\sigma|$, then any daughter site~$\sigma'$
of~$\sigma$ will have its value updated to
\begin{equation}
X_{\sigma'}(t+1)=X_{\sigma'}+\frac1bQ_\sigma^{(\theta)}
=Q_{\sigma'}^{(\theta)}.
\end{equation}
Hence, if the site~$\sigma\in\T_b$ avalanches at time~$t=|\sigma|$, then~$Q_\sigma^{(\theta)}=X_\sigma(t)\ge1$. It
follows that~$\eusmA^{(v)}$, with~$v=\theta-X_\varnothing$, is
the set of sites that are open and connected to the root by a path
of open~sites.
\end{proofsect}

\begin{remark}
\label{rem2} Let us indicate what makes the case~$x_\star=1-\frac1b$ so subtle. Given a sequence~$(c_k)$ of
positive numbers, let us call~$\sigma\in\T_b$ open if~$X_\sigma\ge
x_\star-c_{|\sigma|}b^{-|\sigma|}$ and closed otherwise. Letting
$p_k=\BbbP(X\ge x_\star-c_kb^{-k})$ and supposing, e.g.,~$bp_k=1+k^{-1/2}$, a general result of Lyons~\cite{Lyons} implies
that the open sites percolate. An easy argument shows that if~$\sigma$ is connected to~$\varnothing$ by a path of open sites,
then $Q_\sigma^{(\theta)} \ge1+b^{-k}(v-1-\sum_{\ell\le k}c_\ell)$
for~$v=\theta-X_\varnothing$. Thus, if~$v>1+\sum_{k\ge0}c_k$,
then, by Proposition~\ref{prop2.3}, there is an infinite avalanche
with a non-zero probability.
\end{remark}

On a similar basis, we can write down the necessary and sufficient
conditions for divergence of the expected size of avalanches. The
criterion will be based on the asymptotic growth of the quantity
\begin{equation}
\label{Zntheta}
Z_n(\theta)=\BbbP\bigl(Q^{(\theta)}_k\ge1,\,k=0,\dots,n\bigr),
\qquad n\ge0.
\end{equation}
Notice that~$Z_n(\theta)=0$ whenever~$\theta<1$.

\vbox{
\begin{theorem}
\label{thm2.4} (1) For all~$\theta\ge1$, the limit
\begin{equation}
\label{Zlimit} \zz=\zz(\rho)=\lim_{n\to\infty}Z_n(\theta)^{1/n}
\end{equation}
exists and is independent of~$\theta$.
\newline
(2) For all~$\rho,\rho'\in\MM$, the function~$\alpha\mapsto\zz(\alpha\rho+(1-\alpha)\rho')$ is continuous in~$\alpha\in[0,1]$.
\newline
(3) Let~$\rho\in\MM$ and let~$x_\star$ correspond to~$\rho$ via
\eqref{xstar}. Define~$\zzc=\frac1b$. If~$\zz(\rho)<\zzc$, then $\E_\rho(|\eusmA^{(v)}|)<\infty$ for all~$v\in(0,\infty)$, while
if~$\zz(\rho) >\zzc$, then~$\E_\rho(|\eusmA^{(v)}|)=\infty$ for
all~$v>1-x_\star$.
\end{theorem}
}

Theorem~\ref{thm2.4} defines a free-energy like functional~$\zz$
and gives the characterization of the divergence of~$\chi^{(v)}$,
as already discussed in Section~\ref{sec1.3}. The continuity
statement in part~(2) indicates that the sets of ``avalanching''
and ``non-avalanching'' measures~$\rho\in\MM$ are separated by a
``surface'' (i.e., set of codimension one) of phase transitions.
We will not try to make the latter more precise; our main reason
for including part~(2) is to have an interpretation of the limit~$\zz(\rho)\to\zzc$, which will be needed in the discussion of the
critical behavior. Under additional mild restrictions on~$\rho$,
it will be shown in Section~\ref{sec4} that~$\E_\rho(|\eusmA^{(v)}|)=\infty$ even for the critical measures~$\rho$, i.e.,
those satisfying~$\zz(\rho)=\zzc$.

\begin{proofsect}{Proof of Theorem~\ref{thm2.4}(1)}
We will start with the cases~$\theta=1$ and~$\theta\ge\theta_b$
which are amenable to subadditive-type arguments. Examining~$Z_{n+m}(\theta)$, we may write (by conditioning on~$X_1,\dots,X_m$)
\begin{equation}
\label{ZZrel} Z_{n+m}(\theta)=\E\Bigl(Z_n(Q_m^{(\theta)})
\prod_{j=0}^m\1_{\{Q^{(\theta)}_j\ge1\}}\Bigr).
\end{equation}
Since~$\theta\mapsto Q_n^{(\theta)}$ is manifestly non-decreasing in~$\theta$, so is the event on the right-hand side of
\eqref{Zntheta} and also~$Z_n(\theta)$ itself. Notice that if~$\theta\ge\theta_b$, then~$Q_k^{(\theta)}\le\theta$ for any~$k\ge0$, while if~$\theta=1$, then the conditions in
\eqref{Zntheta} force~$Q_k^{(\theta)}\ge1$. Thus, for~$\theta=1$
we obtain the submultiplicative bound
\begin{equation}
\label{ZZ1sub} Z_{n+m}(1)\ge Z_n(1)Z_m(1),
\end{equation}
while for any~$\theta\ge\theta_b$ we get the supermultiplicative
bound
\begin{equation}
\label{ZZthetasub} Z_{n+m}(\theta)\le Z_n(\theta)Z_m(\theta).
\end{equation}
By standard theorems,~$Z_n(1)^{1/n}$ tends to a limit,~$\zz_1$,
while~$Z_n(\theta)^{1/n}$ for~$\theta\ge\theta_b$ tends to a
(possibly~$\theta$-dependent) limit~$\zz_{\theta}$. Moreover,
\eqref{ZZrel} in fact implies that~$Z_{n+m}(\theta)\le
Z_n(\theta_b+\theta b^{-m})$ and~$\zz_\theta$ is thus constant for
all~$\theta>\theta_b$. We will use~$\zz_\star$ to denote the
common value of~$\zz_\theta$ for~$\theta>\theta_b$. Note that~$Z_n(1)^{1/n}\le\zz_1$ while~$Z_n(\theta)^{1/n}\ge\zz_\star$ for
all~$n\ge1$ and all~$\theta>\theta_b$.

Since~$\theta\mapsto Z_n(\theta)$ is non-decreasing, to prove
\eqref{Zlimit}, we just need to show that~$\zz_\star$ equals~$\zz_1$. If~$x_\star<1-\frac1b$, then~$\zz_\star=0$ and there is
nothing to prove, so let us suppose that~$x_\star\ge1-\frac1b$ for
the rest of the proof. 
As it turns out, we will have to address a number of distinct cases.
These are determined by whether the inequality in~$x_\star\ge1-\frac1b$
is strict or not and by whether the quantity
\begin{equation}
\kappa_\epsilon=\rho\bigl([x_\star-\epsilon,x_\star]\bigr)
\end{equation}
is strictly less than~$\zz_1$ or not for some (particular)~$\epsilon>0$.
Specifically, we will distinguish the following cases:
\settowidth{\leftmargini}{(CASE3333)}
\begin{enumerate}
\item[CASE 1:]
$x_\star>1-\frac1b$ and~$\kappa_\epsilon<\zz_1$ for some~$\epsilon>0$ with~$x_\star-\epsilon > 1 - \frac1b$.
\item[CASE 2:]
$x_\star>1-\frac1b$
but~$\kappa_\epsilon=\zz_1$ for all~$\epsilon>0$ with~$x_\star-\epsilon > 1 - \frac1b$.
\item[CASE 3:]
$x_\star=1-\frac1b$.
\end{enumerate}
As is easily observed, CASE~2 represents the situation where~$\rho$ assigns
no mass to the interval $(1-\frac1b,x_\star)$, while CASE~3 corresponds to the
similar situation when this interval itself is empty.
In view of the trivial inequality~$\zz_1\ge\rho([1-\frac1b,x_\star])$, we must eventually have~$\kappa_\epsilon<\zz_1$ whenever~$\rho$ has any mass in~$(1-\frac1b,x_\star)$.
Hence, the first situation is clearly generic.

In order to address the first two cases (with~$x_\star>1-\frac1b$) we need to establish an
inequality between~$Z_n(\theta)$ and~$\zz_1^n$ for all~$\theta\in[1,\theta_b)$. 
Explicitly, we claim that for~$x_\star>1-\frac1b$ and any~$\theta\in[1,\theta_b)$, 
there is an~$H(\theta)<\infty$ such~that
\begin{equation}
\label{Z1upper} Z_n(\theta) \leq H(\theta)\zz_1^n,\qquad n\ge1.
\end{equation}
Indeed, let~$\epsilon>0$ be such that~$\theta_b-\theta>\epsilon\frac b{b-1}$ and~$x_\star-\epsilon\ge1-\frac1b$ and pick~$m$ so~that
\begin{equation}
\label{2.14a} (x_\star - \epsilon)\Bigl[1 + \frac 1b + \dots +
\frac 1{b^{m-1}}\Bigr] + \frac 1{b^m} \geq \theta.
\end{equation}
Consider the formula \eqref{Zntheta} for~$Z_{n+m}(1)$ but with the
first~$m$ coordinates restricted to the event~$\EE=\{X_1,\dots,X_m\geq x_\star- \epsilon\}$. Notice that on~$\EE$, the conditions involving~$Q_1^{(1)},\dots,Q_m^{(1)}$ are
automatically satisfied. By a derivation similar to
\twoeqref{ZZrel}{ZZ1sub} we have
\begin{equation}
\label{kappabd} Z_{n+m}(1) \geq \kappa_\epsilon^m Z_n(\theta).
\end{equation}
Along with the upper bound~$Z_{n+m}(1) \leq \zz_1^{n+m}$, this
implies \eqref{Z1upper} with~$H(\theta)=(\zz_1/\kappa_\epsilon)^m$.
(We note that, since~$x_\star$ is the suppremum of the support of~$\rho$,
we have~$\kappa_\epsilon>0$ for all~$\epsilon>0$.)

\smallskip
Now we are ready to prove that~$\zz_\star=\zz_1$ in all of the three cases above:

\smallskip\noindent
\textsl{CASE~1}:
Suppose that~$x_\star>1-\frac1b$ and~$\kappa_\epsilon<\zz_1$ for some~$\epsilon>0$ with~$x_\star-\epsilon > 1 - \frac1b$. Let~$\theta>\theta_b$ be small enough that $\theta_\epsilon=x_\star-\epsilon+\frac\theta b<\theta_b$. Then
\begin{equation}
Z_n(\theta) \leq \kappa_\epsilon Z_{n-1}(\theta) +
(1-\kappa_\epsilon)Z_{n-1}(\theta_\epsilon) =
Z_{n-1}(\theta)\biggl[\kappa_\epsilon + (1-\kappa_\epsilon)
\frac{Z_{n-1}(\theta_\epsilon)}{Z_{n-1}(\theta)}\biggr].
\end{equation}
Using \eqref{Z1upper} and the bound~$Z_{n-1}(\theta) \geq
\zz_\star^{n-1}$ we obtain
\begin{equation}
\label{Zthetabig} Z_n(\theta) \leq
Z_{n-1}(\theta)\biggl[\kappa_\epsilon +
(1-\kappa_\epsilon)H(\theta_\epsilon)\Bigl(
\frac{\zz_1}{\zz_\star}\Bigr)^{n-1}\biggr].
\end{equation}
Let~$\kappa_\epsilon(n)$ denote the quantity in the square
brackets, and let us set~$n=2m$ in \eqref{Zthetabig} and iterate
the bound~$m$~times. This gives~$Z_{2m}(\theta) \leq
\kappa_\epsilon(m)^m Z_{m}(\theta)$. If we still entertain the
possibility that~$\zz_1<\zz_\star$, then the~$m\to\infty$ limit
gives~$\zz_\star\le\lim_{m\to\infty}\kappa_\epsilon(m)=\kappa_\epsilon$,
which contradicts the bound~$\zz_1\ge\kappa_\epsilon$. Therefore,
once~$x_\star>1-\frac1b$ and~$\kappa_\epsilon<\zz_1$ for some~$\epsilon>0$
we have~$\zz_1=\zz_\star$.

\medskip\noindent
\textsl{CASE~2}:
Suppose now that~$x_\star>1-\frac1b$
but~$\kappa_\epsilon=\zz_1$ for all~$\epsilon>0$ with~$x_\star-\epsilon > 1 - \frac1b$. 
Notice that this in fact implies that~$\zz_1=\rho(\{x_\star\})$.
We first observe, using
\eqref{kappabd} with~$n=1$, that~$Z_{m+1}(1)\ge\kappa_\epsilon^m\rho([1-\frac1b\theta,x_\star])$
whenever~$\theta$,~$\epsilon$ and~$m$ satisfy \eqref{2.14a}. As a
consequence of \eqref{ZZ1sub}, we have
\begin{equation}
\label{2.19b}
\frac{\zz_1}{\kappa_\epsilon}\ge\biggl[\frac{\rho([1-\frac1b\theta,x_\star])}
{\kappa_\epsilon}\biggr]^{\frac1{m+1}}.
\end{equation}
Now if~$\rho((1-\frac1b\theta_b,x_\star))>0$, we would have $\rho((1-\frac1b\theta_b,x_\star])>\rho(\{x_\star\})=\kappa_\epsilon$ which would by \eqref{2.19b} imply that~$\zz_1>\kappa_\epsilon$, a contradiction. (This fact will be
important later, so we restate it as a corollary right after this
proof.) Hence, we must have~$\rho((1-\frac1b\theta_b,x_\star))=0$. To prove that~$\zz_1=\zz_\star$, let~$\theta>\theta_b$ be small enough that $\theta_0=1+\frac1b(\theta-\theta_b)<\theta_b$. Now either~$X_k=x_\star$ for all~$k=1,\dots,n$, or there is a~$k$ such
that~$X_k\le1-\frac1b\theta_b$. Noting that then~$Q_k^{(\theta)}\le\theta_0$, we thus~have
\begin{equation}
Z_n(\theta)\le \rho\bigl(\{x_\star\}\bigr)^n
+\sum_{k=1}^n\rho\bigl(\{x_\star\}\bigr)^{k-1}
\rho\bigl([0,1-\tfrac1b\theta_b]\bigr)Z_{n-k}(\theta_0).
\end{equation}
Using \eqref{Z1upper}, this gives~$Z_n(\theta)\le
\zz_1^n+n\zz_1^{n-1}\rho([0,1-\frac1b\theta_b])H(\theta_0)$, proving that~$\zz_\star=\zz_1$ holds for~CASE~2 as well.

\medskip\noindent
\textsl{CASE~3}:
Suppose finally that~$x_\star=1-\frac1b$ and note that then~$\theta_b=1$. Immediately, we
have~$Z_n(1)=\rho(\{x_\star\})^n$ and therefore~$\zz_1=\rho(\{x_\star\})$, while for any~$\theta>\theta_b$ we have $\bigcap_{k=0}^n\{Q_k^{(\theta)}\ge1\}\subset\bigcap_{k=1}^n
\{X_k\ge x_\star-b^{-k}(\theta-\theta_b)\}$. Therefore,
\begin{equation}
Z_n(\theta)\le\prod_{k=1}^n\BbbP\bigl(X_k\ge
x_\star-b^{-k}(\theta-\theta_b)\bigr),
\end{equation}
which implies that~$\zz_\star\le\lim_{k\to\infty}\BbbP(X_1\ge
x_\star-b^{-k}(\theta-\theta_b))=\rho(\{x_\star\})=\zz_1$.
\end{proofsect}

This completes the proof of part~(1) of Theorem~\ref{thm2.4}.
As mentioned earlier, we would like to underscore one aspect of the above proof.

\begin{corollary}
\label{cor2.5} Let~$\rho\in\MM$ and suppose
that~$x_\star>1-\frac1b$ and~$\rho((1-\frac1b\theta_b,x_\star))>0$. Then
there is an~$\epsilon>0$ with~$x_\star-\epsilon>\frac{b-1}b$ such that~$\zz(\rho)>\rho([x_\star-\epsilon,x_\star])$.
\end{corollary}

\begin{proofsect}{Proof}
See the argument following \eqref{2.19b}.
\end{proofsect}

Next we will prove the continuity of~$\alpha\mapsto\zz(\alpha\rho+(1-\alpha)\rho')$
as stated in Theorem~\ref{thm2.4}(2):

\begin{proofsect}{Proof of Theorem~\ref{thm2.4}(2)}
Throughout this proof we will write~$Z_n^{(\rho)}(\theta)$ instead
of just~$Z_n(\theta)$ to emphasize the dependence on the
underlying measure~$\rho$. Let~$\rho_0,\rho_1\in\MM$ and let~$\rho_\alpha=(1-\alpha)\rho_0+\alpha\rho_1$. Clearly, to prove
(2), it suffices to show that~$\alpha\mapsto\zz(\rho_\alpha)$ is
right continuous at~$\alpha=0$.

Fix~$\alpha>0$ and let~$(T_k)$ be a sequence of~$0,1$-valued
i.i.d.\ random variables with $\Prob(T_k=0)=\alpha$. Let~$(X_k)$
and~$(X_k')$ be two independent sequences of i.i.d.\ random
variables, both independent of~$(T_k)$, with distributions~$\rho_0^\N$ and~$\rho_1^\N$, respectively. Let~$(X_k^{(\alpha)})$
be the sequence defined~by
\begin{equation}
\label{XXXYrel} 
X_k^{(\alpha)}=T_kX_k+(1-T_k)X_k',\qquad k\ge1.
\end{equation}
Clearly,~$(X_k^{(\alpha)})$ are i.i.d.\ with joint distribution~$\rho_\alpha^\N$. Let us use~$\BbbP_\alpha$ to denote the joint
distribution of~$(X_k)$,~$(X_k')$, and~$(T_k)$.

Let~$Q_{n}^{(\theta,\alpha)}$ be given by \eqref{Qntheta} with~$X_1,\dots,X_n$ replaced by~$X_1^{(\alpha)},\dots,X_n^{(\alpha)}$. Then~$Z_n^{(\rho_\alpha)}(\theta)$ is given by
\eqref{Zntheta} with~$Q_n^{(\theta)}$ replaced by~$Q_{n}^{(\theta,\alpha)}$ and~$\BbbP$ replaced by~$\BbbP_\alpha$.
As will be seen shortly, the main object of interest is the
conditional expectation given the values~$(T_k)$:
\begin{equation}
\label{Znalphatheta}
Z_{n,\alpha}\bigl(\,\theta\,|(T_k)\bigr)=\BbbP_\alpha
\bigl(Q_{\ell}^{(\theta,\alpha)}\ge1,\,\ell=0,\dots,n\big|(T_k)\bigr).
\end{equation}
Indeed, let~$\theta\in[1,\theta_b]$ and, given~$(T_k)$, let~$(I_i)$ be the connected blocks of sites~$k\in\{0,\dots,1\}$ such
that~$T_k=1$ and let~$(J_j)$ be the connected sets of sites with~$T_k=0$. By \eqref{XXXYrel}, the~$X_k^{(\alpha)}$ for~$k\in I_i$
are distributed according to~$\rho_0$, while those for~$k\in J_j$
are distributed according to~$\rho_1$. Then an analogue of
\eqref{ZZrel} for the quantity in \eqref{Znalphatheta} along with
the bounds~$Z_n(1)\le Z_n(\theta)\le Z_n(\theta_b)$ for~$\theta\in[1,\theta_b]$ allow us to conclude that
\begin{equation}
\label{2.16} \prod_i Z_{|I_i|}^{(\rho_0)}(1) \prod_j
Z_{|J_j|}^{(\rho_1)}(1)\le
Z_{n,\alpha}\bigl(\,\theta\,|(T_k)\bigr)\le \prod_i
Z_{|I_i|}^{(\rho_0)}(\theta_b) \prod_j
Z_{|J_j|}^{(\rho_1)}(\theta_b).
\end{equation}

In order to estimate the right hand side of \eqref{2.16}, note
that the existence of the limit in \eqref{Zlimit} implies that for
all~$\delta>0$ there is~$C_\delta\in[1,\infty)$, such that for
both~$\rho=\rho_0$ and~$\rho=\rho_1$,
\begin{equation}
\label{dalsi} Z_n^{(\rho)}(\theta_b)\le C_\delta
(1+\delta)^n\zz(\rho)^n,\qquad n\ge1.
\end{equation}
Let~$\E_\alpha$ denote the expectation with respect to~$\BbbP_\alpha$. Using \eqref{dalsi} in \eqref{2.16}, observing
that the total number of occurrences of~$C_\delta$ is less than~$2k_1(Y)$, where~$k_1(Y)=\sum_j|J_j|$, and noting that~$k_1(Y)$
has the binomial distribution with parameter~$\alpha$ under~$\BbbP_\alpha$ allows us to write
\begin{equation}
\label{2.17} Z_n^{(\rho_\alpha)}(\theta) =\E_\alpha
Z_{n,\alpha}\bigl(\,\theta\,|(T_k)\bigr)\le (1+\delta)^n
\bigl((1-\alpha)\zz(\rho_0)+\alpha C_\delta^2\zz(\rho_1)\bigr)^n.
\end{equation}
By taking~$n\to\infty$, we get~$\lim_{\alpha\downarrow0}\zz(\rho_\alpha)\le
(1+\delta)\zz(\rho_0)$. But~$\delta$ was arbitrary, hence,~$\lim_{\alpha\downarrow0}\zz(\rho_\alpha)\le\zz(\rho_0)$. The
argument for the lower bound,~$\lim_{\alpha\downarrow0}\zz(\rho_\alpha)\ge\zz(\rho_0)$, is
completely analogous.
\end{proofsect}

Finally, we also need to prove part (3) of Theorem~\ref{thm2.4}:

\begin{proofsect}{Proof of Theorem~\ref{thm2.4}(3)}
By Proposition~\ref{prop2.3},~$\sigma\in\eusmA^{(v)}$ is exactly
the event that the path between (and including)~$\sigma$ and~$\varnothing$ consists of sites~$\sigma'$ with~$Q^{(\theta)}_{\sigma'}\ge1$, where~$\theta=X_\varnothing+v$. But
then
\begin{equation}
\BbbP_\rho(\sigma\in\eusmA^{(v)})
=\E_\rho\bigl(Z_{|\sigma|}(X_\varnothing+v)\bigr),
\end{equation}
where the final average is over~$X_\varnothing$. To get the
expected size of~$\eusmA^{(v)}$, we sum over all~$\sigma$,
\begin{equation}
\label{2.20} \E_\rho\bigl(|\eusmA^{(v)}|\bigr)
=\sum_{n\ge0}\E_\rho\bigl(Z_n(X_\varnothing+v)\bigr)b^n.
\end{equation}
The existence of the limit~$Z_n^{1/n}(\theta)$ independent of~$\theta$ (for~$\theta\ge1$) tells us that~$\E_\rho(|\eusmA^{(v)}|)<\infty$ whenever~$\zz(\rho)<\zzc$, while~$\E_\rho(|\eusmA^{(v)}|)=\infty$ once~$\zz(\rho)>\zzc$ and~$v>1-x_\star$.
\end{proofsect}

\section{Absence of intermediate phase}
\label{sec3}
\subsection{Sharpness of phase transition}
\label{sec3.1}\noindent The goal of this section is to show that
the phase transitions defined by presence/absence of an infinite
avalanche and divergence of avalanche size occur at the same
``point,''~$\zzc=\frac1b$. This rules out the possibility of an
intermediate phase. Moreover, we will prove that the transition is
\textit{second order} in the sense that there is no infinite avalanche at~$\zz=\zzc$.

Unfortunately, our proof will require certain restrictions on the underlying measure~$\rho$.
The delicate portion of~$[0,1]$ is the region~$I=[1-\frac1b\theta_b,\frac{b-1}b)$.
Clearly, some conditions are needed to ensure  that there is not too much mass at the right-end of~$I$---i.e., that $\rho([(1-\frac1b-\epsilon,1-\frac1b))\to0$ sufficiently fast as $\epsilon\downarrow0$---to avoid the sort of counterexamples described in
Remarks~\ref{rem1} and~\ref{rem2} of Section~\ref{sec2}.
With a lot of additional work than what is about to hit, all of the forthcoming can be proved under the assumption that~$\rho$ has an~$L^p$ density, for a~$p>1$, in the interval~$I$.
However, this requires dealing with ``singularities'' in the region above~$I$.
(The region below~$I$ is of no consequence because any directed path in the avalanched set can only harbor a finite number of values from this set.)
Notwithstanding, most of the interesting mathematics---with only a fraction of unpleasant technicalities---is captured by assuming that the measure~$\rho$ has a bounded density.

\begin{definition}
Let~$\MMsharp$ be the set of Borel probability measures~$\rho$ on~$[0,1]$ that are absolutely continuous with respect to Lebesgue measure
on~$[0,1]$ with the associated density~$\phi_\rho$ bounded in~$L^\infty$ norm on~$[0,1]$, i.e.,~$\Vert\phi_\rho\Vert_\infty<\infty$, and that obey $\rho([1-\frac1b,1])>0$.
\end{definition}

The requirement that~$\rho$ has no positive mass in~$[1-\frac1b,1]$ represents no additional loss of generality since the opposite case, namely $x_\star<1-\frac1b$, actually has~$\zz(\rho)=0$ and is therefore far away from having an avalanche (see Theorem~\ref{thm2.4}).
It is worth noting that~$\MMsharp$ is a convex subset of~$\MM$.
The ability to take convex
combinations of elements of~$\MMsharp$ will be crucial in the
discussion of the critical behavior, see Section~\ref{sec4}.

\smallskip
Our second main theorem is then as follows:

\begin{theorem}
\label{thm3.1} Suppose that~$\rho\in\MMsharp$ and define~$\zzc=\frac1b$.
\newline\indent
(1) If~$\zz(\rho)\le\zzc$, then~$\BbbP_\rho(|\eusmA^{(v)}|=\infty)=0$ for all~$v\in(0,\infty)$.
\newline\indent
(2) If~$\zz(\rho)>\zzc$, then~$\BbbP_\rho(|\eusmA^{(v)}|=\infty)>0$ for all~$v\in(1-x_\star,\infty)$.
\end{theorem}

The proof of Theorem~\ref{thm3.1} requires introducing two
auxiliary random variables~$V_\infty$ and~$Q_\infty$. These will
be defined in next two subsections, the proof is therefore
deferred to Section~\ref{sec3.4}. The random variable~$Q_\infty$
will be a cornerstone of our analysis of the critical process, see
Section~\ref{sec4}. The underlying significance of both~$V_\infty$
and~$Q_\infty$ is the distributional identity that each of them
satisfies.

\subsection{Definition of~$V_\infty$}
\label{sec3.2}\noindent In this section we define a random
variable~$V_\infty$ which is, roughly speaking, the minimal value
of~$v$ that needs to be added to the root in order to trigger an
infinite avalanche. For~$n\ge1$ let
\begin{equation}
\label{Vn} V_n=\inf\bigl\{v\in(0,\infty)\colon
\boldX^{(v)}(t+1)\ne\boldX^{(v)}(t),\, t=0,\dots,n-1\bigr\}.
\end{equation}
(A logical extension of this definition to~$n=0$ is~$V_0\equiv
0$.) In plain words, if~$v\ge V_n$, then the avalanche process
will propagate to at least the~$n$-th level. Clearly,~$V_n$ is an
increasing sequence; we let~$V_\infty$ denote the~$n\to\infty$
limit of~$V_n$. Formally,~$V_\infty$ can be infinite; in fact,
since the event~$\{V_\infty<\infty\}$ is clearly a tail event,~$\BbbP_\rho(V_\infty<\infty)$ is either one or zero.

Let us use~$\Psi_n$ to denote the distribution function of~$V_n$,
i.e.,
\begin{equation}
\label{Psindef} \Psi_n(\vartheta)=\BbbP_\rho(V_n\le\vartheta).
\end{equation}
The aforementioned properties of~$V_n$ lead us to a few immediate
observations about~$\Psi_n$: First,~$\Psi_n$ is a decreasing
sequence of non-decreasing functions. Second, the limit
\begin{equation}
\label{Psidef} \Psi(\vartheta)=\lim_{n\to\infty}\Psi_n(\vartheta),
\end{equation}
exists for all~$\vartheta\in(0,\infty)$ and~$\Psi(\vartheta)=\BbbP_\rho(V_\infty\le\vartheta)$. Third, we have~$\Psi\not\equiv0$ if and only if~$\BbbP_\rho(V_\infty<\infty)=1$.
Moreover, each of~$\Psi_n$ is in principle computable:

\begin{lemma}
\label{lemma3.2} 
Let $\rho\in\MM$.
Then the sequence~$(\Psi_n)$ satisfies the recurrence equation
\begin{equation}
\label{Psirecur}
\Psi_{n+1}(\vartheta)=\E_\rho\biggl(\Phi_b\Bigl(\Psi_n\bigl(
\textstyle\frac{X_\varnothing+\vartheta}b\bigr)\displaystyle\Bigr)
\1_{\{X_\varnothing\ge1-\vartheta\}}\biggr), \qquad n\ge0,
\end{equation}
where~$\Psi_0(\vartheta)=\1_{\{\vartheta\ge0\}}$ and
\begin{equation}
\label{Phib} \Phi_b(y)=1-(1-y)^b, \qquad 0\le y\le1.
\end{equation}
\end{lemma}

\begin{proofsect}{Proof}
Let~$\T_b^{(\sigma)}$ denote the subtree of~$\T_b$ rooted at~$\sigma$ and let~$V_n^{(\sigma)}$ denote the random variable
defined in the same way as~$V_n$ but here for the tree~$\T_b^{(\sigma)}$. Then we have
\begin{equation}
\label{vnrecur}
\{V_{n+1}\le\vartheta\}=\{X_\varnothing\ge1-\vartheta\}\cap\Bigl\{
\,\min_{\sigma\in\{1,\dots,b\}}V_n^{(\sigma)}\le\frac{X_\varnothing+\vartheta}b
\Bigr\}.
\end{equation}
But for all~$\sigma\in\{1,\dots,b\}$, the~$V_n^{(\sigma)}$'s are
i.i.d.\ with common distribution function~$\Psi_n$, so we have
\begin{equation}
\label{Phirel} \BbbP_\rho\bigl(\,\min_{\sigma\in\{1,\dots,b\}}
V_n^{(\sigma)}\le\vartheta\bigr)
=\Phi_b\bigl(\Psi_n(\vartheta)\bigr).
\end{equation}
From here the claim follows by noting that~$V_n^{(\sigma)}$ are
independent of~$X_\varnothing$.
\end{proofsect}

\begin{corollary}
\label{cor3.3} 
Let $\rho\in\MM$.
Then the distribution function of~$V_\infty$ satisfies
the equation
\begin{equation}
\label{Psieq} \Psi(\vartheta)=\E_\rho\biggl(\Phi_b\Bigl(\Psi\bigl(
\textstyle\frac{X_\varnothing+\vartheta}b\bigr)\displaystyle\Bigr)
\1_{\{X_\varnothing\ge1-\vartheta\}}\biggr), \qquad\vartheta\ge0.
\end{equation}
\end{corollary}
\begin{proofsect}{Proof}
This is an easy consequence of \eqref{Psirecur} and the Bounded
Convergence Theorem.
\end{proofsect}

On the basis of \eqref{Psieq} and some percolation arguments, the
answer to the important question whether~$\Psi\equiv0$ or not can
be given by checking whether~$\Psi(\vartheta)=0$ for reasonable
values of~$\vartheta$:

\begin{proposition}
\label{prop3.4} 
Let $\rho\in\MMsharp$. 
Suppose that~$\Psi\not\equiv0$. Then
\begin{equation}
\label{essrel}
\inf\bigl\{\vartheta\ge0\colon\Psi(\vartheta)>0\bigr\}=1-x_\star.
\end{equation}
\end{proposition}

\begin{proofsect}{Proof}
Let~$\vartheta_\star$ denote the infimum on the left-hand side of
\eqref{essrel}.
Note that $x_\star>1-\frac1b$ by $\rho\in\MMsharp$.
Since~$\rho$ is absolutely continuous with respect to the Lebesgue 
measure on~$[0,1]$, there is an $\eta>0$ such that $x_\star-\eta>1-\frac1b$ and~$\rho([x_\star-\eta,x_\star])<\frac1b$. Now~$\frac1b$ is the threshold for the site percolation on~$\T_b$, so the sites with~$X_\sigma>x_\star-\eta$ do not percolate. Let~$\G_n=\{\sigma\in\T_b\colon|\sigma|=n\}$ be the~$n$-th generation of~$\T_b$. Pick two integers~$N,N'$ such that~$N'\ge N$ and let~$\HH_{N,N'}$ be the event that~$\G_N$ and~$\G_{N'}$ are separated by a ``barrier'' of sites~$\sigma$ with~$X_\sigma\le x_\star-\eta$. By taking~$N'\gg N\gg 1$, the probability of~$\HH_{N,N'}$ can be made as close to one as desired.

Let~$\vartheta>\vartheta_\star$ and pick~$N_0$ so large that~$\vartheta b^{-N_0}$ is less than~$\frac\eta2$. Find~$N,N'\ge N_0$ such  that
$1-\BbbP_\rho(\HH_{N,N'})$ is strictly smaller than~$\BbbP_\rho(|\eusmA^{(\vartheta)}|=\infty)$, i.e., we have
$\BbbP_\rho(\{|\eusmA^{(\vartheta)}|=\infty\}\cap\HH_{N,N'})>0$.
Now for any~$\epsilon\in(0,\frac\eta2)$, we will produce a configuration with an infinite avalanche that has a starting value~$v=1-x_\star+\epsilon$. Draw a configuration~$(\bar X_\sigma)$ subject to the constraint that~$\bar X_\sigma\ge x_\star-\epsilon$ for all~$\sigma\in\T_b$ with~$|\sigma|\le N'$. Let~$(X_\sigma)$ belong to the set~$\{|\eusmA^{(\vartheta)}|=\infty\}\cap\HH_{N,N'}$ and define~$X_\sigma'$ by putting
\begin{equation}
X_\sigma'=\begin{cases}
\bar X_\sigma\vee X_\sigma,\qquad&\text{if }|\sigma|\le N',
\\
X_\sigma,\qquad&\text{otherwise}.
\end{cases}
\end{equation}
Let~$X_\sigma^{\prime,(v)}(t)$ denote the process corresponding to the initial configuration~$(X_\sigma')$ and initial value $v>0$, and let~$X^{(\vartheta)}_\sigma(t)$ be the corresponding process for~$(X_\sigma)$ and~$\vartheta$. Let~$\eusmA^{\prime,(v)}$ and~$\eusmA^{(\vartheta)}$ be the corresponding avalanche sets.

The configuration~$(X_\sigma)$ exhibits an infinite avalanche, so there is a site~$\sigma$ on one of the aforementioned ``barriers'' separating~$\G_N$ and~$\G_{N'}$, which belongs to an infinite oriented path inside~$\eusmA^{(\vartheta)}$. By the assumption that~$x_\star-\eta>1-\frac1b$ it is clear that, if~$v>1-x_\star+\epsilon$ and~$t=|\sigma|$, then~$\eusmA^{\prime,(v)}$ will reach~$\sigma$. But~$X_{\sigma'}'\ge X_{\sigma}$ for all sites on the path from~$\varnothing$ to~$\sigma$, so we have 
\begin{equation}
\label{pred}
X_\sigma^{\prime,(v)}(t)-X^{(\vartheta)}_\sigma(t)
\ge\eta-\epsilon-\frac{\vartheta-v}{b^N}>0,
\end{equation}
where we used that~$b^N\vartheta\le\frac\eta2$ and~$\epsilon<\frac\eta2$ to derive the last inequality.
Now the set~$\eusmA^{(\vartheta)}$ contains a path from~$\sigma$ to infinity and, by \eqref{pred} and~$X_{\sigma'}'\ge X_{\sigma'}$ 
for~$\sigma'$ ``beyond''~$\sigma$, 
this path will also be contained in~$\eusmA^{\prime,(v)}$. 
Consequently, an infinite avalanche will occur in configuration~$(X_\sigma')$ starting from a value~$v>1-x_\star+\epsilon$ whenever it did in configuration~$(X_\sigma)$ starting from~$\vartheta$. This establishes~$\vartheta_\star=1-x_\star$, as claimed.
\end{proofsect}

\subsection{Definition of~$Q_\infty$}
\label{sec3.3}\noindent The second random variable, denoted by~$Q_\infty$ is a limiting version of the objects~$Q_n^{(\theta)}$
defined in \eqref{Qntheta}. Let~$Y=(Y_1,Y_2,\dots)$ be a sequence
of i.i.d.\ random variables with joint distribution $\BbbP=\rho^\N$. These are, in a certain sense, the same
quantities as the~$X$'s discussed earlier, however, the~$Y$'s will be
ordered in the opposite way. Similarly to \eqref{Qntheta}, let
\begin{equation}
\label{Qntheta1}
Q_{n,k}^{(\theta)}=Y_k+\frac1bY_{k+1}+\dots+\frac1{b^{n-k}}
Y_n+\frac\theta{b^{n-k+1}}, \qquad 1\le k\le n.
\end{equation}
For completeness, we also let~$Q_{0,1}^{(\theta)}=\theta$. 

Let~$\BB$ be the Borel~$\sigma$-algebra on~$[0,1]^\N$ equipped
with the standard product topology. Suppose that~$\rho([\frac{b-1}b,1])>0$---which is assured if $\rho\in\MMsharp$. For any~$n\ge1$ and~$\theta\ge1$, let~$\BbbP_n^{(\theta)}$ be the conditional law on~$\BB$ defined by
\begin{equation}
\label{Pntheta}
\BbbP_n^{(\theta)}(\,\cdot\,)=\BbbP\bigl(\,\cdot\,\big|\,
Q^{(\theta)}_{n,\ell}\ge1,\,\ell=2,\dots,n\bigr),
\end{equation}
The latter is well defined because~$\BbbP(Q^{(\theta)}_{n,\ell}\ge1)>0$ for all~$\ell=2,\dots,n$,~$\{Q^{(\theta)}_{n,\ell}\ge1\}$ are increasing and~$\BbbP(\cdot)$
is FKG. Intentionally, the variable~$Y_1$ is not constrained by
the conditioning in \eqref{Pntheta}.

\smallskip
Next we give conditions for the existence of the limiting law~$\lim_{n\to\infty}\BbbP_n^{(\theta)}$:

\begin{proposition}
\label{prop3.5} 
Let~$\rho\in\MMsharp$ and let~$\theta_0>\theta_b$. 
Then there exists numbers~$A=A(\rho,\theta_0)<\infty$ and~$\zeta=\zeta(\rho)>0$ such
that for all bounded measurable functions~$f=f(Y_1,\dots,Y_k)$
and all $\theta,\theta'\in[1,\theta_0]$,
\begin{equation}
\label{cbd}
\bigl|\E^{(\theta)}_{n+1}(f)-\E^{(\theta')}_n(f)\bigr|\le A
e^{-\zeta(n-k)}\Vert f\Vert_\infty,\qquad n\ge k.
\end{equation}
In particular, whenever~$\theta\ge1$, the limit law
\begin{equation}
\label{limeq}
\wBbbP(\cdot)=\lim_{n\to\infty}\BbbP^{(\theta)}_n(\cdot)
\end{equation}
exists and is independent of~$\theta$. 
Moreover, the quantities~$A(\rho,\theta_0)$ and~$\zeta(\rho)$ are bounded away from infinity and zero uniformly in any convex set~$\NN\subset\MMsharp$ with finitely many extreme points.
\end{proposition}

The proof of Proposition~\ref{prop3.5} uses a coupling argument,
which requires some rather extensive preparations and is therefore
deferred to Section~\ref{sec5}. (The actual proof appears at the
end of Section~\ref{sec5.3}.)

\smallskip
We will use~$\wE$ to denote the expectation with respect to~$\wBbbP$ whenever the latter is well defined. Let us define a
random variable~$Q_\infty$ on~$([0,1]^\N,\BB,\wBbbP)$ by the
formula
\begin{equation}
\label{Qinfty} Q_\infty=\sum_{k\ge1}\frac{Y_k}{b^{k-1}}.
\end{equation}
Notice that~$Q_\infty$ is supported in~$[\frac1b,\theta_b]$,
because~$Y_1$ is not constrained by the conditioning in
\eqref{Pntheta}. 

\begin{corollary}
\label{Holder-cor}
Let~$\rho\in\MMsharp$ and let~$\theta\ge1$. Let~$Q_{n,1}^{(\theta)}$ be as in~\eqref{Qntheta1}, where the
variables $Y_1,\dots,Y_n$ are distributed according
to~$\BbbP_n^{(\theta)}$. Then~$Q_{n,1}^{(\theta)}$ tends to~$Q_\infty$ in distribution 
as~$n\to\infty$.
Moreover, for each~$\theta_0>\theta_b$ 
and each~$C<\infty$
there are constants~$D=D(\rho,\theta_0)<\infty$ and~$\varsigma=\varsigma(\rho)>0$ such that if
$f(\theta)$ is a function obeying the Lipschitz bound on~$[0,\theta_0]$,
\begin{equation}
\label{HBBbd}
\bigl|f(\theta)-f(\theta')\bigr|\le C\Vert f\Vert_\infty\,|\theta-\theta'|,
\qquad \theta,\theta'\in[0,\theta_0],
\end{equation}
where~$\Vert f\Vert_\infty=\sup_{\theta\le\theta_0}|f(\theta)|$, 
then
\begin{equation}
\label{Hbd}
\Bigl|\,\E_{n}^{(\theta)}\bigl(f(Q_{n,1}^{(\theta)})\bigr)
-\wE\bigl(f(Q_\infty)\bigr)\Bigr|\le D\,
\Vert f\Vert_\infty e^{-\varsigma n}
\end{equation}
holds for all~$\theta\in[1,\theta_0]$.
The quantities~$D(\rho,\theta_0)$
and~$\varsigma(\rho)$ are bounded away from infinity and zero uniformly in any convex set~$\NN\subset\MMsharp$ with finitely many extreme points.
\end{corollary}

The proof of Corollary~\ref{Holder-cor} is given in Section~\ref{sec5.4}.
As already mentioned, a principal tool for our later investigations will
be the distributional identity for~$Q_\infty$ stated
below.

\begin{proposition}
\label{prop3.6} Let~$\rho\in\MMsharp$. If~$X$ is a random variable with law~$\BbbP=\rho$, independent of~$Q_\infty$,~then
\begin{equation}
\label{disteq} \BbbP\otimes\wBbbP
\biggl(X+\frac{Q_\infty}b\in\cdot\,\bigg|\,
Q_\infty\ge1\biggr)=\wBbbP(Q_\infty\in\cdot\,).
\end{equation}
\end{proposition}

The proof of Proposition~\ref{prop3.6} will also be given in
Section~\ref{sec5}. Proposition~\ref{prop3.5} and the proof of
Proposition~\ref{prop3.6} immediately yield an extension of
Theorem~\ref{thm2.4}(1), stated as Corollary~\ref{cor3.7}, which
will also be useful in subsequent developments. The proof of
Corollary~\ref{cor3.7} is given in Section~\ref{sec5.4}.

\begin{corollary}
\label{cor3.7} 
Suppose that~$\rho\in\MMsharp$. Then~$\zz(\rho)=\wBbbP(Q_\infty\ge1)$. Moreover, the limit
\begin{equation}
\label{psirho}
\psi_\rho(\theta)=\lim_{n\to\infty}Z_n(\theta)\zz(\rho)^{-n}
\end{equation}
exists for all~$\theta\ge0$ and, 
for all~$\theta_0>\theta_b$, 
there are~$A'=A'(\rho,\theta_0)<\infty$ and~$\zeta'=\zeta'(\rho)>0$ such~that
\begin{equation}
\label{error}
\bigl|Z_n(\theta)\zz(\rho)^{-n}-\psi_\rho(\theta)\bigr|\le A' e^{-\zeta' n}
\end{equation}
holds for all~$\theta\in[0,\theta_0]$ and all~$n\ge1$.
Furthermore, the function~$\psi_\rho$
has the following properties:
\begin{enumerate}
\item[(1)]
$\psi_\rho(\theta)\in(0,\infty)$
for all~$\theta\ge1$ while~$
\psi_\rho(\theta)=0$ for~$\theta<1$.
\item[(2)]
$\theta\mapsto\psi_\rho(\theta)$
is non-decreasing and Lipschitz continuous for all~$\theta\ge1$. More precisely, 
there is a~$C=C(\rho,\theta_0)<\infty$ such that
$|\psi_\rho(\theta)-\psi_\rho(\theta')|\le 
C\psi_\rho(\theta_0)|\theta-\theta'|$ for all~$\theta,\theta'\in[1,\theta_0]$.
\item[(3)]
If~$\rho,\rho'\in\MMsharp$ and~$\rho_\alpha=(1-\alpha)\rho+\alpha\rho'$ for each~$\alpha\in[0,1]$, then~$\alpha\mapsto\psi_{\rho_\alpha}(\theta)$
is continuous in~$\alpha\in[0,1]$ for~all~$\theta\ge0$.
\end{enumerate}
The quantities~$A'(\rho,\theta_0)$,~$\zeta'(\rho)$ and~$C(\rho,\theta_0)$ are bounded away from infinity and zero uniformly in any convex set~$\NN\subset\MMsharp$ with finitely many extreme points.
\end{corollary}

\begin{remark}
The Lipschitz continuity of~$\theta\mapsto\psi_\rho(\theta)$ is a direct consequence of our assumption that~$\rho$ has a \textit{bounded} density~$\phi_\rho$ with respect to the Lebesgue measure on~$[0,1]$. If~$\phi_\rho$ is only in~$L^p([0,1])$ for some~$p>1$, then the appropriate concept will be H\"older continuity with a~$p$-dependent exponent. The same will be true for various other Lipschitz continuous quantities later in this paper.
\end{remark}

\subsection{Proof of Theorem~\ref{thm3.1}}
\label{sec3.4}\noindent With random variable~$Q_\infty$ at our
disposal, the sharpness of the phase transition in our avalanche
model is almost immediate.

\begin{proofsect}{Proof of Theorem~\ref{thm3.1}}
Let~$\rho\in\MMsharp$ and abbreviate~$\zz=\zz(\rho)$. Let~$x_\star$ be as in \eqref{xstar}. We begin by introducing the
quantity
\begin{equation}
G_n=\wE\Bigl(\Psi_n\bigl(\textstyle\frac{Q_\infty}b\bigr)
\displaystyle\1_{\{Q_\infty\ge1\}}\Bigr).
\end{equation}
The recursive equation \eqref{Psirecur} and
Proposition~\ref{prop3.6} then give
\begin{equation}
\label{AAeq}
\begin{aligned}
G_{n+1}&=\wBbbP(Q_\infty\ge1)\,\E_\rho\otimes\wE
\biggl(\Psi_n\Bigl(\frac{\scriptstyle
X_\varnothing+\frac1b{Q_\infty}}b\Bigr)
\1_{\{X_\varnothing+\frac1b{Q_\infty}\ge1\}}\,\Big|\,Q_\infty\ge1\biggr)\\
&=\zz\,\wE\Bigl(\Phi_b\Bigl(
\Psi_n\bigl(\textstyle\frac{Q_\infty}b
\bigr)\1_{\{Q_\infty\ge1\}}\Bigr)\Bigr),
\end{aligned}
\end{equation}
where we have used the fact that~$\zz=\wBbbP(Q_\infty\ge1)$ from
Corollary~\ref{cor3.7}.

Let us first analyze the cases~$b\zz\le1$. By using Jensen's
inequality in \eqref{AAeq} we get that
\begin{equation}
\label{AAupper} G_{n+1}\le\zz\,\Phi_b(G_n)\le\frac1b\Phi_b(G_n).
\end{equation}
An inspection of the graph of~$y\mapsto\Phi_b(y)$ reveals that if
\eqref{AAupper} holds, then~$G_n\to0$. By
Proposition~\ref{prop3.4}, this is compatible with~$\Psi\not\equiv0$ only if
\begin{equation}
\label{aa} \frac{Q_\infty}b\le1-x_\star\qquad\wBbbP\text{-almost
surely}.
\end{equation}
However, a simple argument shows that~$\esssup Q_\infty=
x_\star\frac b{b-1}$ whenever~$\zz>0$. This contradicts
\eqref{aa}, because~$x_\star>\frac {b-1}b$ (as implied by~$\rho\in\MMsharp$) forces~$1-x_\star<x_\star\frac b{b-1}$. Thus, if~$b\zz\le 1$, then~$\Psi$
must be identically zero.

Next we will attend to the cases~$b\zz>1$. We will suppose that~$\Psi_n\to0$ and work to derive a contradiction. Since~$n\mapsto\Psi_n$ is a monotone sequence of monotone functions, the
convergence to~$\Psi$ is uniform on~$[0,1]$ and, in particular, on
the range of values that~$\frac1bQ_\infty$ takes. Using that~$\Phi_b(y)\ge by-\frac12b(b-1)y^2$ for all~$y\in[0,1]$ and
invoking \eqref{AAeq}, we can write
\begin{equation}
\label{AAlower} G_{n+1}\ge b\zz(1-\epsilon_n)G_n,
\end{equation}
where~$\epsilon_n=\frac12(b-1)\Psi_n(1)$. Since~$b\zz>1$ and~$\epsilon_n\to0$, we have~$G_{n+1}\ge G_n$ for~$n$ large enough.
An inspection of \eqref{Psirecur} shows that, since~$x_\star>1-\frac1b$, we have~$\Psi_n(\vartheta)>0$ for all~$\vartheta>\frac1b$. Hence~$G_n>0$ for all~$n\ge0$. But then
\eqref{AAlower} forces~$G_n$ to stay uniformly bounded away from
zero, in contradiction with our assumption that~$G_n\to 0$.
Therefore, once~$b\zz>1$, we must have~$\Psi\not\equiv0$.
\end{proofsect}

\vbox{
\section{Critical behavior}
\label{sec4}\nopagebreak
\subsection{Critical exponents}
\label{sec4.1}\noindent In this section we establish, under
certain conditions on~$\rho$, the essential behavior of the model
at the critical point~$\zzc=\frac1b$. In particular, we describe
the asymptotics for the critical distribution of avalanche sizes,
the power law behavior for the probability of an infinite
avalanche as~$\zz\downarrow\zzc$ and, finally, the exponent for
the divergence of~$\chi^{(v)}$ as~$\zz\uparrow\zzc$.
}

\begin{theorem}
\label{thm4.1} Let~$\rho\in\MMsharp$ and let~$x_\star$ be as in
\eqref{xstar}. Suppose~$\zz(\rho)=\zzc$, where~$\zzc=\frac1b$.
Then there are functions~$\tau,\TT\colon(1-x_\star,\infty)\to(0,\infty)$ and~$\varTheta\colon[0,\infty)\to[0,\infty)$ such that the following
holds:
\newline
(1) If~$\rho'\in\MMsharp$ and~$\rho_\alpha=\alpha\rho'+(1-\alpha)\rho$ satisfies~$\zz(\rho_\alpha)<\zzc$ for all~$\alpha\in(0,1]$, then for all~$v>1-x_\star$,
\begin{equation}
\E_{\rho_\alpha}\bigl(|\eusmA^{(v)}|\bigr)=\frac{\tau(v)}{\zzc-\zz(\rho_\alpha)}
\bigl[1+o(1)\bigr],\qquad \alpha\downarrow0.
\end{equation}
\newline
(2) For all~$v\ge0$,
\begin{equation}
\label{critscal} \BbbP_\rho\bigl(|\eusmA^{(v)}|\ge
n\bigr)=\frac{\varTheta(v)}{n^{1/2}}\bigl[1+o(1)\bigr], \qquad
n\to\infty,
\end{equation}
where~$\varTheta(v)>0$ for~$v>1-x_\star$.
\newline
(3) If~$\rho'\in\MMsharp$ and~$\rho_\alpha=\alpha\rho'+(1-\alpha)\rho$ satisfies~$\zz(\rho_\alpha)>\zzc$ for all~$\alpha\in(0,1]$, then for all~$v>1-x_\star$,
\begin{equation}
\label{superscal} \BbbP_{\rho_\alpha} \bigl(|\eusmA^{(v)}|=\infty\bigr) =\bigl(\zz(\rho_\alpha)-\zzc\bigr)
\,\TT(v)\bigl[1+o(1)\bigr],\qquad \alpha\downarrow0.
\end{equation}
\end{theorem}

\begin{remark}
The proof of Theorem~\ref{thm4.1} makes frequent use of the properties of the random variable~$Q_\infty$ defined in Section~\ref{sec3.2}. The relevant statements are Propositions~\ref{prop3.5} and~\ref{prop3.6} and Corollaries~\ref{Holder-cor} and~\ref{cor3.7}, whose proofs come only in Section~\ref{sec5}. Modulo these claims, Section~\ref{sec4} is essentially self-contained and can be read without a reference to Section~\ref{sec5}.
\end{remark}

Part~(1) of Theorem~\ref{thm4.1} can be proved based on the
already-available information; the other parts will require some
preparations and their proofs are postponed to the next section.

\begin{proofsect}{Proof of Theorem~\ref{thm4.1}(1)}
Let~$\rho,\rho'\in\MMsharp$ be such that~$\zz(\rho_\alpha)<\zzc=\zz(\rho)$ for~$\rho_\alpha=(1-\alpha)\rho+\alpha\rho'$  and all~$\alpha\in(0,1]$. Let~$\chi^{(v)}(\alpha)=\E_{\rho_\alpha}(|\eusmA^{(v)}|)$. By \eqref{2.20},
\begin{equation}
\label{4.4}
\chi^{(v)}(\alpha)=\sum_{n\ge0}\E_{\rho_\alpha}\bigl(
Z_n^{(\rho_\alpha)}(X_\varnothing+v)\bigr)\,b^n,
\end{equation}
where~$\E_{\rho_\alpha}$ is the expectation with respect to~$X_\varnothing$ in~$\rho_\alpha$ and~$Z_n^{(\rho_\alpha)}$ is
defined by \eqref{Zntheta} using~$\rho_\alpha$.

In order to estimate the sum we will use~$A''$ and~$\zeta''$ to denote the
worst case scenarios for the quantities~$A'(\rho_\alpha,\theta_0)$ and~$\zeta'(\rho_\alpha)$ from Corollary~\ref{cor3.7}. Explicitly, we let $A''=\sup_{0\le\alpha\le1}A'(\rho_\alpha,\theta_0)$ and~$\zeta''=\inf_{0\le\alpha\le1}\zeta'(\rho_\alpha)$, where $\theta_0>\theta_b$ is to be determined shortly. 
Note that $A''<\infty$ and $\zeta''>0$ by uniformity of the bounds on $A'(\rho_\alpha,\theta_0)$ and $\zeta'(\rho_\alpha)$ in the convex set $\NN=\{\rho_\alpha\colon\alpha\in[0,1]\}$.
Then we have, for all~$n\ge1$ and
all~$\theta\in[1,\theta_0]$,
\begin{equation}
\label{4.5}
b^nZ_n^{(\rho_\alpha)}(\theta)
=b^n\zz(\rho_\alpha)^n\psi_{\rho_\alpha}(\theta)+b^n\zz(\rho_\alpha)^n E_n(\theta),
\end{equation}
where~$\psi_{\rho_\alpha}(\theta)$ is as in \eqref{psirho} while~$E_n(\theta)$ is the ``error term.'' 
Using the bounds from Corollary~\ref{cor3.7},~$E_n(\theta)$ is estimated by
$|E_n(\theta)|\le A''e^{-\zeta'' n}$.
By continuity of~$\alpha\mapsto\psi_{\rho_\alpha}(\theta)$, we get
\begin{equation}
\label{4.6}
\sum_{n\ge0}b^n Z_n^{(\rho_\alpha)}(\theta)=
\frac{\psi_\rho(\theta)+o(1)}{1-b\zz(\rho_\alpha)},
\end{equation}
where~$o(1)$ tends to zero as~$\alpha\downarrow0$ uniformly on compact sets of~$\theta\in[1,\theta_0]$.

Let~$\tau(v)=b^{-1}\E_\rho(\psi_\rho(X_\varnothing+v))$ and note that~$\tau(v)>0$ for all~$v>1-x_\star$. 
Let us take the maximum of $x_\star+v$ and $2\theta_b$ for the quantity~$\theta_0$ above.
Then \eqref{4.4} and \eqref{4.6} imply
\begin{equation}
\label{4.7} \chi^{(v)}(\alpha)=\frac{\tau(v)}{\zzc-\zz(\rho_\alpha)}
\bigl[1+o(1)\bigr],
\end{equation}
where~$o(1)$ tends to zero as~$\alpha\downarrow0$, for all~$v\ge1-x_\star$. 
\end{proofsect}

It remains to establish parts (2) and (3) of Theorem~\ref{thm4.1}.
To ease derivations, instead of looking at the asymptotic size of~$\eusmA^{(v)}$, we will focus on a slightly different set:
\begin{equation}
\eusmB^{(\theta)}=\begin{cases} \{\varnothing\},\qquad&\text{if
}\eusmA^{(\theta-X_\varnothing)}=\emptyset,
\\
\bigl\{\sigma\in\T_b\colon m(\sigma)\in\eusmA^{(\theta-X_\varnothing)}\bigr\},\qquad&\text{otherwise.}
\end{cases}
\end{equation}
(Here we take~$\eusmA^{(\theta')}=\emptyset$ whenever~$\theta'<1$.) Clearly,~$\eusmB^{(\theta)}$ is the original
avalanche set together with its boundary (i.e., the set of sites in~$\T_b$, where the avalanche has ``spilled'' some material). Since
both sets are connected and both contain the root (with the
exception of the case~$\eusmA^{(\theta-X_\varnothing)}=\emptyset$), their sizes satisfy the
relation:
\begin{equation}
\label{BandA} |\eusmB^{(\theta)}|=(b-1)|\eusmA^{(\theta-X_\varnothing)}|+1.
\end{equation}
(This relation holds even if~$\eusmA^{(\theta-X_\varnothing)}=\emptyset$.) The asymptotic
probability of the events~$\{|\eusmA^{(v)}|\ge n\}$ as~$n\to\infty$ is thus basically equivalent to that of~$\{|\eusmB^{(\theta)}|\ge(b-1)n\}$.

\subsection{Avalanches in an external field}
\label{sec4.2}\noindent 
Following a route which is often used in the analysis of critical systems,
our proof of Theorem~\ref{thm4.1} will be accomplished by the addition of
extra degrees of freedom that play the role of an \textit{external field}.
Let~$\lambda\in[0,1]$ be fixed and let us color each
site of~$\T_b$ ``green'' with probability~$\lambda$. Given a
realization of this process, let~$\eusmG$ denote the random set
of ``green'' sites in~$\T_b$. Let~$\BbbP_{\rho,\lambda}(\cdot)$ be
the joint probability distribution of the ``green'' sites and~$(X_\sigma)$. The principal quantity of interest is~then
\begin{equation}
B_\infty(\theta,\lambda)=\BbbP_{\rho,\lambda} \bigl(\eusmB^{(\theta)}\cap \eusmG\ne\emptyset\bigr).
\end{equation}
It is easy to check that, as~$\lambda\downarrow0$, the number~$B_\infty(\theta,\lambda)$ tends to the probability~$\BbbP_\rho(|\eusmB^{(\theta)}|=\infty)$. In particular,
Theorem~\ref{thm3.1} guarantees that~$B_\infty(\theta,\lambda)\to0$ as~$\lambda\downarrow0$ if~$\zz(\rho)\le\zzc$, while~$B_\infty(\theta,\lambda)$ stays
uniformly positive as~$\lambda\downarrow0$ when~$\zz(\rho)>\zzc$
and~$\theta\ge1$.

Let~$\psi_\rho(\theta)$ be as in Corollary~\ref{cor3.7} and let~$c_\rho\in(0,\infty)$ be the quantity defined by
\begin{equation}
\label{crho} \frac1{c_\rho^2}=\frac{b-1}2\,\wE\Bigl(\bigl[\E
\bigl(\psi_\rho
\bigl(X\!+\!\textstyle\frac{Q_\infty}b\bigr)\bigr)\bigr]^2\,\Big|\,
Q_\infty\ge1\Bigr).
\end{equation}
Here~$X$ and~$Q_\infty$ are independent with distributions~$\BbbP=\rho$ and~$\wBbbP$, respectively. It turns out that the
asymptotics of~$B_\infty(\theta,\lambda)$ for critical~$\rho$ can
be described very precisely:

\begin{proposition}
\label{prop4.2} Let~$\rho\in\MMsharp$ satisfy~$\zz(\rho)=\zzc$.
For each~$\theta\in(0,\infty)$,
\begin{equation}
\label{4.11a}
\lim_{\lambda\downarrow0}\frac{B_\infty(\theta,\lambda)}
{\sqrt\lambda}=c_\rho\psi_\rho(\theta).
\end{equation}
\end{proposition}

\smallskip
Proposition~\ref{prop4.2} is proved in Section~\ref{sec4.4}. Now
we are ready to prove Theorem~\ref{thm4.1}(2):

\begin{proofsect}{Proof of Theorem~\ref{thm4.1}(2)}
We begin by noting the identity
\begin{equation}
\frac{B_\infty(\theta,\lambda)}\lambda
=\sum_{n\ge1}\BbbP_\rho\bigl(|\eusmB^{(\theta)}|\ge
n\bigr)(1-\lambda)^{n-1},\qquad \lambda\in(0,1],
\end{equation}
which is derived by expressing $\BbbP_\rho(|\eusmB^{(\theta)}|=n)$ as the difference between $\BbbP_\rho(|\eusmB^{(\theta)}|=n)$ and $\BbbP_\rho(|\eusmB^{(\theta)}|=n+1)$.
Since~$B_\infty(\theta,\lambda)=\sqrt\lambda\,
(c_\rho\psi_\rho(\theta)+o(1))$ as~$\lambda\downarrow0$ and since~$n\mapsto \BbbP_\rho(|\eusmB^{(\theta)}|\ge n)$ is a decreasing
sequence, standard Tauberian theorems (e.g., Karamata's Tauberian
Theorem for Power Series, see Corollary~1.7.3 in~\cite{BGT})
guarantee that
\begin{equation}
\label{Blim} \BbbP_\rho(|\eusmB^{(\theta)}|\ge n)=
c_\rho\frac{\psi_\rho(\theta)}{\Gamma(\tfrac12)}\frac1{\sqrt
n}\bigl[1+o(1)\bigr],\qquad n\to\infty,
\end{equation}
(Strictly speaking, the above Tauberian theorem applies only when~$\psi_\rho(\theta)>0$; in the opposite case, i.e., when~$\theta<1$, we have~$\eusmB^{(\theta)}=\{\varnothing\}$ and 
there is nothing to prove.) In order to
obtain the corresponding asymptotics for~$\BbbP_\rho(|\eusmA^{(v)}|\ge n)$, we first note that, by \eqref{BandA},
\begin{equation}
\BbbP_\rho\bigl(|\eusmA^{(v)}|\ge n\bigr)=\BbbP_\rho\bigl(|\eusmB^{(X_\varnothing+v)}|\ge (b-1)n+1\bigr).
\end{equation}
By applying \eqref{Blim} on the right-hand side and invoking the
Bounded Convergence Theorem, we immediately get the desired
formula \eqref{critscal} with
\begin{equation}
\varTheta(v)= \frac{c_\rho}{(b-1)^{1/2}\Gamma(\tfrac12)}\,
\E_\rho\bigl(\psi_\rho(X_\varnothing+v)\bigr),
\end{equation}
where~$\E_\rho$ is the expectation over~$X_\varnothing$. Clearly,~$v\mapsto\varTheta(v)$ is non-decreasing because~$\theta\mapsto\psi_\rho(\theta)$ is non-decreasing, while~$\varTheta(v)>0$ for~$v>1-x_\star$ because~$\psi_\rho(\theta)>0$
for~$\theta\ge1$.
\end{proofsect}

Similarly we can also describe the asymptotics of~$\BbbP_\rho(|\eusmB^{(\theta)}|=\infty)$ as~$\zz(\rho)\downarrow\zzc$:

\vbox{
\begin{proposition}
\label{prop4.3} Let~$\rho,\rho'\in\MMsharp$ and define~$\rho_\alpha=(1-\alpha)\rho+\alpha\rho'$. Suppose that~$\zz(\rho)=\zzc$ and~$\zz(\rho_\alpha)>\zzc$ for all~$\alpha\in(0,1]$. Then for all~$\theta\in(0,\infty)$,
\begin{equation}
\label{4.11b} \frac{\BbbP_\rho(|\eusmB^{(\theta)}|=\infty)}
{\zz(\rho_\alpha)-\zzc}=bc_\rho^2\psi_\rho(\theta)+o(1), \qquad
\alpha\downarrow0,
\end{equation}
where~$\psi_\rho(\theta)$ is as in Corollary~\ref{cor3.7} and~$c_\rho$
is as in \eqref{crho}.
\end{proposition}
}

Proposition~\ref{prop4.3} is proved in Section~\ref{sec4.5}. Now
we are ready finish the proof of Theorem~\ref{thm4.1}(3):

\begin{proofsect}{Proof of Theorem~\ref{thm4.1}(3)}
By \eqref{BandA} we clearly have that
\begin{equation}
\BbbP_\rho\bigl(|\eusmA^{(v)}|=\infty\bigr)=\BbbP_\rho\bigl(|\eusmB^{(X_\varnothing+v)}|=\infty\bigr).
\end{equation}
By conditioning on~$X_\varnothing+v=\theta$ and invoking
\eqref{4.11b}, we can easily derive that the asymptotic formula
\eqref{superscal} holds with~$\TT$ given by~$\TT(v)=bc_\rho^2\E_\rho (\psi_\rho(X_\varnothing+v))$.
\end{proofsect}

As we have seen, Propositions~\ref{prop4.2} and~\ref{prop4.3} have
been instrumental in the proof of Theorem~\ref{thm4.1}(2) and~(3).
The following three sections are devoted to the proofs of the two
propositions. After some preliminary estimates, which constitute a
substantial part of Section~\ref{sec4.3}, we will proceed to
establish the critical asymptotics (Section~\ref{sec4.4}). The
supercritical cases can then be handled along very much the same
lines of argument, the necessary changes are listed in
Section~\ref{sec4.5}.

\subsection{Preliminaries}
\label{sec4.3}\noindent This section collects some facts about the
quantity~$B_\infty(\theta,\lambda)$ and its~$\theta$ and~$\lambda$
dependence. We begin by proving a simple identity for~$B_\infty(\theta,\lambda)$:

\begin{lemma}
\label{lemma4.2} Let~$\rho\in\MM$ and let~$\Phi_b$ be as in
\eqref{Phib}. Then
\begin{equation}
\label{LB}
B_\infty(\theta,\lambda)=\lambda+(1-\lambda)\1_{\{\theta\ge1\}}\,\Phi_b\bigl(\E_\rho
B_\infty\bigl(X_\varnothing+\textstyle\frac1b \displaystyle
\theta,\lambda \bigr)\bigr),
\end{equation}
\end{lemma}

\begin{proofsect}{Proof}
If~$\theta<1$, then~$B_\infty(\theta,\lambda)=\lambda$ and
\eqref{LB} clearly holds true. Let us therefore suppose that~$\theta\ge1$. Let~$\eusmB_\sigma^{(\theta)}$ denote the object~$\eusmB^{(\theta)}$ for the subtree of~$\T_b$ rooted at~$\sigma$.
Then
\begin{equation}
\label{LBset} \bigl\{\eusmB^{(\theta)}\cap \eusmG\ne\emptyset\bigr\} =\{\varnothing\in\eusmG\}
\cup\biggl(\{\varnothing\not\in \eusmG\}\cap\bigcup_{\sigma=1}^b
\bigl\{\eusmB_\sigma^{(X_\sigma+\frac1b\theta)}\cap \eusmG\ne\emptyset\bigr\}\biggr).
\end{equation}
The claim then follows by using the independence of the sets in
the large parentheses on the right hand side of \eqref{LBset}
under the measure~$\BbbP_{\rho,\lambda}(\cdot)$.
\end{proofsect}

Our next claim concerns continuity properties of~$B_\infty(\theta,\lambda)$
as a function of~$\theta$:

\begin{lemma}
\label{lemma4.3} 
For each~$\rho\in\MMsharp$ satisfying~$\zz(\rho)<\zzc e$ and each~$\theta_0>\theta_b$ there is a~$C=C(\rho,\theta_0)<\infty$ such that
\begin{equation}
\label{HHbd}
\bigl|B_\infty(\theta,\lambda)-B_\infty(\theta',\lambda)\bigr|\le
C  B_\infty(\theta_0,\lambda)|\theta-\theta'|
\end{equation}
for all~$\lambda\ge0$ and all~$\theta,\theta'\in[1,\theta_0]$. 
The bound~$C(\rho,\theta_0)<\infty$ is uniform
in any convex set~$\NN\subset\{\rho\in\MMsharp\colon \zz(\rho)<\zzc 
e\}$ with finitely many extreme~points.
\end{lemma}

\begin{proofsect}{Proof}
Let us assume that~$\theta\ge\theta'$. To derive \eqref{HHbd},
we will regard~$B_\infty(\theta,\lambda)$ and~$B_\infty(\theta',\lambda)$ as
originating from the same realization of~$(X_\sigma)$ and the
``green'' sites. Then~$\Delta=B_\infty(\theta,\lambda)-B_\infty(\theta',\lambda)$ is
dominated by the probability (under~$\BbbP_{\rho,\lambda}$) that
there is a site~$\sigma\in\T_b$,~$\sigma\ne\varnothing$,
with the properties:
\begin{enumerate}
\item[(1)]~$Q_{\sigma'}^{(\theta')}\ge1$ for all~$\sigma'=m^k(\sigma)$ with~$k=1,\dots,|\sigma|$.
\item[(2)]~$Q_\sigma^{(\theta')}<1$ but~$Q_\sigma^{(\theta)}\ge1$.
\item[(3)]~$\eusmB_\sigma^{(\theta_0)}\cap\eusmG\ne\emptyset$, where~$\eusmB_\sigma^{(\theta_0)}$ is the  set~$\eusmB^{(\theta_0)}$ for
the subtree~$\T_b^{(\sigma)}$ rooted at~$\sigma$.
\end{enumerate}
Indeed, any realization of~$(X_\sigma)$ and the ``green'' sites
contributing to~$\Delta$ obeys~$\eusmB^{(\theta')}\cap
\eusmG=\emptyset$ and~$\eusmB^{(\theta)}\cap\eusmG\ne\emptyset$. But then
there must be a site~$\sigma$ on the inner boundary of~$\eusmB^{(\theta')}$ where the avalanche corresponding to~$\theta'$
stops but that corresponding to~$\theta$ goes on. (Since~$\theta,\theta'\ge1$, we must have~$\sigma\ne\varnothing$.) Consequently,~$Q_{\sigma'}^{(\theta')}\ge1$ for any~$\sigma'$ on the path
connecting~$\sigma$ to the root, but~$Q_\sigma^{(\theta')}<1\le
Q_\sigma^{(\theta)}$, justifying conditions~(1) and~(2) above.
Since~$Q_\sigma^{(\theta)}\le\theta_0$, and since the~$\theta$-avalanche continuing on from~$\sigma$ must eventually
reach a ``green'' site, we see that also condition~(3) above must
hold.

Let~$\rho\in\MMsharp$ be such that~$\zz(\rho)<\zzc e$. 
Using the independence of the events
described in (1), (2) and~(3) above, and recalling the definitions \eqref{Zntheta} and \eqref{Pntheta}, we can thus estimate
\begin{equation}
\label{Deltabd} \Delta\le
B_\infty(\theta_0,\lambda)\sum_{\sigma\in\T_b\smallsetminus\{\varnothing\}}
Z_{|\sigma|-1}(\theta')\,
\BbbP_{|\sigma|}^{(\theta')}
\bigl(Q_{|\sigma|,1}^{(\theta)}\ge1>Q_{|\sigma|,1}^{(\theta')}\bigr).
\end{equation}
Abbreviate~$K_n(\theta,\theta')=
\BbbP_n^{(\theta')}(Q_{n,1}^{(\theta)}\ge1>Q_{n,1}^{(\theta')})$. Since~$Y_1$ is independent of all the other~$Y$'s in the measure~$\BbbP_n^{(\theta')}$, we have
\begin{equation}
K_n(\theta,\theta')=
\bigl\{\rho([1-\tfrac{\vartheta'}b,1-\tfrac{\vartheta}b))\colon
\vartheta'-\vartheta\le|\theta-\theta'|b^{-n+1}\bigr\}.
\end{equation}
Here~$\vartheta$, resp.,~$\vartheta'$ play the role of~$Q_{n,2}^{(\theta)}$, resp.,~$Q_{n,2}^{(\theta')}$ and the interval in the argument of~$\rho$ exactly corresponds to the inequalities~$Q_{n,1}^{(\theta)}=Y_1+\frac1b\vartheta\ge1>Y_1+\frac1b\vartheta'=Q_{n,1}^{(\theta')}$.

To estimate the supremum, we recall that~$\rho(\textd x)=\phi_\rho(x)\textd x$ where~$\phi_\rho$ is bounded. Then
\begin{equation}
K_n(\theta,\theta')\le\Vert\phi_\rho\Vert_\infty\,
|\theta-\theta'|\, b^{-n+1},
\qquad n\in\N.
\end{equation}
Now, by
Corollary~\ref{cor3.7},~$Z_n(\theta)\le C\zz(\rho)^n$
for some~$C<\infty$ uniformly in~$\rho$ on convex sets~$\NN\subset\{\rho\in\MMsharp\colon \zz(\rho)<\zzc e\}$
with finitely many extreme points and uniformly in~$\theta\le\theta_0$. 
Therefore, the right-hand side of
\eqref{Deltabd} is bounded by~$B_\infty(\theta_0,\lambda)|\theta-\theta'|$ times a
sum that converges whenever~$\zz(\rho)<\zzc e$, uniformly in~$\rho\in\NN$, where~$\NN$ is as above. 
This proves the desired
claim.
\end{proofsect}

Let~$\rho\in\MMsharp$ and let~$Q_\infty$ be the random variable
defined in Section~\ref{sec3.3}, independent of both the green
sites and~$X_\sigma$. Let us introduce the quantity
\begin{equation}
B^\star_\infty(\lambda)
=\wE\bigl(B_\infty(Q_\infty,\lambda)\bigr),
\end{equation}
The significance of~$B^\star_\infty(\lambda)$ is that it
represents a stationary form of~$B_\infty(\cdot,\lambda)$, i.e.,~$B^\star_\infty(\lambda)$ is a very good approximation of the
probability~$\BbbP_{\rho,\lambda}(\eusmB_\sigma^{(\theta')}\cap
\eusmG=\emptyset\,|\, \sigma\in\eusmA^{(v)})$, where~$\theta'=Q_\sigma^{(X_\varnothing+v)}$ and where~$\eusmB_\sigma^{(\theta')}$ is the quantity~$\eusmB^{(\theta')}$ for
trees rooted at~$\sigma$ very far from~$\varnothing$. Let
\begin{equation}
\label{kapparho} \varkappa_\rho(\lambda)
=\wE\Bigl(\bigl[\E\bigl( B_\infty
\bigl(X\!+\!\textstyle\frac{Q_\infty}b,
\lambda\bigr)\bigr)\bigr]^2\,\Big|\, Q_\infty\ge1\Bigr),
\end{equation}
where~$X$ and~$Q_\infty$ are independent with distributions~$\BbbP=\rho$ and~$\wBbbP$, respectively.
For critical distributions,~$B_\infty^\star(\lambda)$ and~$\varkappa_\rho(\lambda)$ are related as follows:

\begin{lemma}
\label{lemma4.3a} Let~$\rho\in\MMsharp$ be such that~$\zz(\rho)=\zzc$. Then
\begin{equation}
\label{LBbd} 
B_\infty^\star(\lambda) =1-\frac{b-1}{2\lambda}\,
\varkappa_\rho(\lambda)\bigl[1+o(1)\bigr],
\qquad\lambda\downarrow0.
\end{equation}
\end{lemma}

\begin{proofsect}{Proof}
Since~$B_\infty(\theta,\lambda)\to0$ as~$\lambda\downarrow0$, we
can expand~$\Phi_b$ on the right hand side of \eqref{LB} to the
second order of Taylor expansion, use that~$\zz(\rho)=\wBbbP(Q_\infty\ge1)$ and apply~$b\zz(\rho)=1$ with the~result
\begin{equation}
\label{4.26} B_\infty^\star(\lambda)=\lambda+
(1-\lambda)B_\infty^\star(\lambda) -\frac{b-1}{2}\,
\varkappa_\rho(\lambda)\bigl[1+o(1)\bigr],
\qquad\lambda\downarrow0.
\end{equation}
(Here we noted that~$B_\infty(X+\frac1bQ_\infty)\le
B_\infty(\theta_b)$ allows us to estimate the error in the Taylor
expansion by~$\varkappa_\rho(\lambda)B_\infty(\theta_b)O(1)$,
which is~$\varkappa_\rho(\lambda)o(1)$ as~$\lambda\downarrow0$.)
Subtracting~$(1-\lambda)B_\infty^\star(\lambda)$ from both sides and
dividing by~$\lambda$, \eqref{LBbd} follows.
\end{proofsect}

Note that, by the resulting expression \eqref{LBbd},~$\varkappa_\rho(\lambda)/\lambda$ tends to a definite limit as~$\lambda\downarrow0$. In the supercritical cases, on the other
hand, Lemma~\ref{lemma4.3a} gets replaced by the following claim:

\begin{lemma}
\label{lemma4.3b} Let~$\rho,\rho'\in\MMsharp$ and define~$\rho_\alpha=(1-\alpha)\rho+\alpha\rho'$. Suppose that~$\zz(\rho)=\zzc$ and~$\zz(\rho_\alpha)>\zzc$ for all~$\alpha\in(0,1]$. Let~$B_\infty^\star(0,\alpha)$ denote the
quantity~$B_\infty^\star(0)$ for the underlying measure~$\rho_\alpha$. Then
\begin{equation}
\label{LBbdsuper}
B_\infty^\star(0,\alpha)=\frac{b-1}{2b}\frac{\varkappa_{\rho_\alpha}(0)}
{\zz(\rho_\alpha)-\zzc}\bigl[1+o(1)\bigr],\qquad
\alpha\downarrow0.
\end{equation}
\end{lemma}

\begin{proofsect}{Proof}
As in Lemma~\ref{lemma4.3a}, we use that~$B_\infty(\theta,0,\alpha)\to0$ as~$\alpha\to0$, where~$B_\infty(\theta,0,\alpha)$ denotes the quantity~$B_\infty(\theta,0)$ for the underlying measure~$\rho_\alpha$.
However, instead of \eqref{4.26}, this time we get
\begin{equation}
\label{4.28}
B_\infty^\star(0,\alpha)\bigl(1-b\zz(\rho_\alpha)\bigr)=
-\frac{b-1}{2}\, \varkappa_{\rho_\alpha}(0)\bigl[1+o(1)\bigr],
\qquad\alpha\downarrow0,
\end{equation}
where we again used that the error in the Taylor approximation can
be bounded by~$\varkappa_{\rho_\alpha}o(1)$. Dividing by~$\zz(\rho_\alpha)-\zzc\ne0$, \eqref{LBbdsuper} follows.
\end{proofsect}

\subsection{Critical asymptotics}
\label{sec4.4}\noindent The purpose of this section is to finally
give the proof of Proposition~\ref{prop4.2}. We begin by proving
an appropriate upper bound on~$B_\infty(\theta,\lambda)$. Note
that, despite being used only marginally, equation \eqref{LBbd} is
a key ingredient of the proof.

\vbox{
\begin{lemma}
\label{lemma4.6} Let~$\rho\in\MMsharp$ satisfy~$\zz(\rho)=\zzc$.
For each~$\theta\ge1$ there is a~$K(\theta)\in(0,\infty)$
such that
\begin{equation}
\label{Abd} \limsup_{\lambda\downarrow0}
\frac{B_\infty(\theta,\lambda)}{\sqrt\lambda}\le K(\theta).
\end{equation}
\end{lemma}
}

\begin{proofsect}{Proof}
Let~$\zz=\zz(\rho)$. We begin by proving \eqref{Abd} for~$\theta=1$. Let
\begin{equation}
\iota(\rho)
=\wE\Bigl(\rho\bigl([1-\frac1bQ_\infty,1])^2\Big|Q_\infty\ge1\Bigr)
\end{equation}
and recall the
definition of~$\varkappa_\rho(\lambda)$ in \eqref{kapparho}. Using
the inequality $B_\infty(\theta,\lambda)\ge
B_\infty(1,\lambda)\1_{\{\theta\ge1\}}$ we derive $\varkappa_\rho(\lambda)\ge \iota(\rho)B_\infty(1,\lambda)^2$.
Inserting this in \eqref{LBbd}, we have
\begin{equation}
B_\infty^\star(\lambda)
\le1- \frac{b-1}{2\lambda}
\iota(\rho)\,B_\infty(1,\lambda)^2\bigl[1+o(1)\bigr],\quad
\lambda\downarrow0.
\end{equation}
Since the left-hand side is always non-negative, \eqref{Abd} for~$\theta=1$ follows with~$K(1)^{-2}=\frac{b-1}2\iota(\rho)$.

Next we will show that for any~$\theta<\theta_b$,~$B_\infty(\theta,\lambda)$ is bounded above by a
($\theta$-dependent) multiple of~$B_\infty(1,\lambda)$. Indeed,
pick an~$\epsilon>0$ such that~$\theta_b-\theta>\epsilon\frac b{b-1}$ and let~$m$ be so large
that \eqref{2.14a} holds. Fix a directed path of~$m$ steps
in~$\T_b$ starting from the root. By conditioning on the event
that~$X_\sigma\ge x_\star-\epsilon$ for all~$\sigma\ne\varnothing$
in the path, we have~$B_\infty(1,\lambda)\ge\rho([x_\star-\epsilon,1])^m
B_\infty(\theta,\lambda)$, i.e.,
\begin{equation}
B_\infty(\theta,\lambda)\le C(\theta)
B_\infty(1,\lambda),\qquad\theta<\theta_b,
\end{equation}
with~$C(\theta)=\rho([x_\star-\epsilon,1])^{-m}<\infty$.

As the third step we prove that \eqref{Abd} holds for values~$\theta$ in slight excess of~$\theta_b$. (The reader will notice
slight similarities with the latter portion of the proof 
of Theorem~\ref{thm2.4}(1).) 
Let $\epsilon>0$ be such
that~$x_\star-\epsilon>1-\frac1b$. By
Corollary~\ref{cor2.5} and the fact that~$\rho\in\MMsharp$, 
we can assume that~$\kappa_\epsilon=\rho([x_\star-\epsilon,x_\star])<\zz$. 
If~$\theta>\theta_b$ is such that~$\theta_\epsilon=x_\star-\epsilon+\frac1b\theta<\theta_b$, then
\eqref{LB} and the bound~$\Phi_b(y)\le by$ imply
\begin{equation}
B_\infty(\theta,\lambda)\le\lambda+(1-\lambda)b\bigl[
\kappa_\epsilon B_\infty(\theta,\lambda)+(1-\kappa_\epsilon)
B_\infty(\theta_\epsilon,\lambda)\Bigr],
\end{equation}
because~$X+\frac1b\theta\le\theta$ for all~$X$ in the support of~$\rho$. Since~$(1-\lambda)b\kappa_\epsilon<b\kappa_\epsilon<1$, we
have
\begin{equation}
B_\infty(\theta,\lambda)\le
\frac{\lambda+(1-\lambda)(1-\kappa_\epsilon)b C(\theta_\epsilon)
B_\infty(1,\lambda)}{1-(1-\lambda)b\kappa_\epsilon}.
\end{equation}
Dividing by~$\sqrt\lambda$ and taking~$\lambda\downarrow0$,
\eqref{Abd} follows with~$K(\theta)=b(1-\kappa_\epsilon)C(\theta_\epsilon)
K(1)/(1-b\kappa_\epsilon)$.

Finally, it remains to prove \eqref{Abd} for general~$\theta\ge\theta_b$. But for that we just need to observe that
\begin{equation}
B_\infty(\theta,\lambda)\le[1-(1-\lambda)^{b^{k+1}}]+
(1-\lambda)^{b^{k+1}}B_\infty(\theta_b+\theta b^{-k},\lambda)
\end{equation}
as follows by conditioning on the first~$k$ layers of~$\T_b$ to be
green-free. By taking~$k$ large enough, $\theta_b+\theta b^{-k}$
is arbitrary close to~$\theta_b$, so the result follows by the
preceding arguments.
\end{proofsect}

Lemma~\ref{lemma4.6} allows us to write the following expression for~$B_\infty(\theta,\lambda)$:

\begin{lemma}
\label{lemma4.9a} Let~$\rho\in\MMsharp$ satisfy~$\zz(\rho)=\zzc$.
Let~$\epsilon(\lambda,\theta)$ be defined by
\begin{equation}
\label{Alimasymp} B_\infty(\theta,\lambda)= \psi_\rho(\theta)
B_\infty^\star(\lambda)+\epsilon(\lambda,\theta),
\end{equation}
where~$\psi_\rho(\theta)$ is as in \eqref{psirho}. Then~$\lim_{\lambda\downarrow0}\epsilon(\lambda,\theta)\lambda^{-1/2}=0$
uniformly on compact sets of~$\theta$.
\end{lemma}

\begin{proofsect}{Proof}
Recall the notation~$Q_{n,1}^{(\theta)}$ from \eqref{Qntheta1},
and let~$\E_n^{(\theta)}$ denote the expectation with respect to
the measure~$\BbbP_n^{(\theta)}$ in \eqref{Pntheta}. We will first
show that
\begin{equation}
\label{Aasymp} B_\infty(\theta,\lambda)= Z_n(\theta)b^n\,
\E_n^{(\theta)}\bigl(B_\infty(Q_{n,1}^{(\theta)},
\lambda)\bigr)+\tilde\epsilon_n(\lambda)
\end{equation}
holds with an~$\tilde\epsilon_n(\lambda)$ satisfying~$\lim_{\lambda\downarrow0}\tilde\epsilon_n(\lambda)\lambda^{-1/2}=0$
for all~$n\ge1$. Let~$\G_n$ denote the~$n$-th generation of~$\T_b$, i.e.,~$\G_n=\{\sigma\in\T_b\colon|\sigma|=n\}$, and let~$\BbbH_n=\bigcup_{m<n}\G_m$.
Recall the notation~$\eusmB_\sigma^{(\theta)}$ for the object~$\eusmB^{(\theta)}$ on the subtree~$\T_b^{(\sigma)}$ of~$\T_b$
rooted at~$\sigma$ and let~$Q_\sigma^{(\theta)}$ be as described
in \eqref{Qsigma}. Given a~$\sigma\in\G_n$, let~$\pi(\sigma)=\{m^k(\sigma)\colon k=1,\dots,n\}$ be the path
of connecting~$\sigma$ to the root.

A moment's thought reveals that, if~$\eusmG\cap\BbbH_n=\emptyset$
(i.e., if there are no green sites in the first~$n-1$ generations
of~$\T_b$), then in order for~$\eusmB^{(\theta)}\cap \eusmG\ne\emptyset$
to occur, the following must hold: First, 
there is a~$\sigma\in\G_n$, such that~$Q_{\sigma'}^{(\theta)}\ge1$ for all~$\sigma'\in\pi(\sigma)$. Second,
the avalanche starting from this~$\sigma$ with an initial
amount~$Q_\sigma^{(\theta)}$ reaches~$\eusmG$. Introducing the
event
\begin{equation}
\label{UUn}
\UU_n=\bigcup_{\sigma\in\G_n}\Bigl(\bigl\{\eusmB_\sigma^{(Q_\sigma^{(\theta)})}\cap \eusmG\ne\emptyset\bigr\}\cap\!\!\bigcap_{\sigma'\in\pi(\sigma)}
\{Q_{\sigma'}^{(\theta)}\ge1\}\Bigr),
\end{equation}
we thus have
\begin{equation}
\label{setrel} 
\BbbP_{\rho,\lambda}(\UU_n)\le B_\infty(\theta,\lambda)\le
\BbbP_{\rho,\lambda}(\UU_n)+\BbbP_{\rho,\lambda}\bigl(
\{\eusmG\cap\BbbH_n\ne\emptyset\}\bigr).
\end{equation}
Since~$\BbbP_\rho(\eusmG\cap\BbbH_n\ne\emptyset)=O(\lambda)$, it
clearly suffices to show that~$\BbbP_{\lambda,\rho}(\UU_n)$ has
the same asymptotics as claimed on the right-hand side of \eqref{Aasymp}.

Since~$\UU_n$ is the union of~$b^n$ events with the same probability, the upper bound 
\begin{equation}
\BbbP_{\rho,\lambda}(\UU_n)\le b^n Z_n(\theta)\E_n^{(\theta)}\bigl(B_\infty(Q_{n,1}^{(\theta)},\lambda)\bigr)
\end{equation}
directly follows using the identity
\begin{equation}
\E_\rho\Bigl(B_\infty(Q_\sigma^{(\theta)},\lambda)
\!\prod_{\sigma'\in\pi(\sigma)}\!
\1_{\{Q_{\sigma'}^{(\theta)}\ge1\}}\Bigr)=
Z_n(\theta)\E_n^{(\theta)}\bigl(B_\infty(Q_{n,1}^{(\theta)},\lambda)\bigr).
\end{equation}
To derive the lower bound, we use the inclusion-exclusion formula. 
The exclusion term (i.e., the sum over intersections of pairs of events from
the union in \eqref{UUn}) is estimated, using the bound in
Lemma~\ref{lemma4.6}, to be less than~$K(\bar\theta)^2b^{2n}\lambda$, where~$\bar\theta=\theta\vee\theta_b$. This proves \eqref{Aasymp}.

Since~$\zz(\rho)b=1$, Corollary~\ref{cor3.7} tells us that~$Z_n(\theta)b^n=\psi_\rho(\theta)+o(1)$. 
The final task is to show that~$\E_n^{(\theta)}(B_\infty(Q_{n,1}^{(\theta)},
\lambda))$ can safely be replaced by its limiting version,~$B_\infty^\star(\lambda)$.
We cannot use Corollary~\ref{Holder-cor} directly, because~$\theta\mapsto B_\infty(\theta,\lambda)$ is known to be Lipschitz continuous only for~$\theta\ge1$. However, by Lemma~\ref{lemma4.2} we know that~$B_\infty(\theta,\lambda)=\lambda$ for~$\theta<1$, which means that we can~write
\begin{equation}
\label{4.45a}
B_\infty(\theta,\lambda)=B_\infty\bigl(\theta\vee1,\lambda)
+\bigl[\lambda-B_\infty(1,\lambda)\bigr]\1_{\{\theta<1\}}.
\end{equation}
Now,~$B_\infty^1(\theta,\lambda)=B_\infty(\theta\vee1,\lambda)$ is Lipschitz continuous in~$\theta$ for all~$\theta\ge0$, so by \eqref{HHbd} and \eqref{Hbd},
\begin{equation}
\label{refbd}
\Bigl|\E_n^{(\theta)}\bigl(B^1_\infty(Q_{n,1}^{(\theta)},
\lambda)\bigr)-\wE\bigl(B^1_\infty(Q_\infty,\lambda)\bigl)\Bigr|\le
DB_\infty(\bar\theta,\lambda)e^{-\varsigma n}
\end{equation}
where~$\varsigma>0$ and~$D=D(\bar\theta)<\infty$. To estimate the contribution of the second term in \eqref{4.45a}, we first note that~$\lambda-B_\infty(1,\lambda)$ is a constant bounded between~$-B_\infty(\bar\theta,\lambda)$ and zero. Hence, we thus need to estimate the difference~$\BbbP_n^{(\theta)}(Q_{n,1}^{(\theta)}<1)-\wBbbP(Q_\infty<1)$. But that can be done using Proposition~\ref{prop3.5}: Let~$k=\lfloor\frac n2\rfloor$ and use the monotonicity of~$\theta\mapsto Q_{k,1}^{(\theta)}$ and \eqref{cbd} to estimate
\begin{equation}
\label{refbd2}
\bigl|\BbbP_n^{(\theta)}(Q_{n,1}^{(\theta)}<1)-\wBbbP(Q_\infty<1)\bigr|
\le
\BbbP_n^{(\theta)}(Q_{k,1}^{(\bar\theta)}\ge 1)-\wBbbP(Q_{k,1}^{(1)}\ge 1)\le
A''e^{-\zeta(n-k)},
\end{equation}
where~$A''=A/(1-e^{-\zeta})$. 
By combining all the previous estimates and invoking \eqref{Abd}, we find that
the difference~$\E_n^{(\theta)}(B_\infty(Q_{n,1}^{(\theta)},
\lambda))-B_\infty^\star(\lambda)$ is proportional to~$e^{-\varsigma'n}\sqrt\lambda$,
where~$\varsigma'>0$.
Using this back in \eqref{Aasymp} 
the claim follows by taking the limits~$\lambda\downarrow0$ and~$n\to\infty$.
\end{proofsect}

Lemmas~\ref{lemma4.6} and~\ref{lemma4.9a} finally allow us to
prove Proposition~\ref{prop4.2}:

\begin{proofsect}{Proof of Proposition~\ref{prop4.2}}
Note that, by using \eqref{Alimasymp} in
\eqref{kapparho} and the definition of~$c_\rho$ in \eqref{crho}, we have
\begin{equation}
\frac{b-1}2\varkappa_\rho(\lambda)=B_\infty^\star(\lambda)^2
c_\rho^{-2}+o(\lambda),\qquad\lambda\downarrow0.
\end{equation}
Then the fact that~$B_\infty^\star(\lambda)$ tends to zero as~$\lambda\downarrow0$ forces, in light of \eqref{LBbd}, that~$\frac{b-1}{2\lambda}\varkappa_\rho(\lambda)\to1$ as~$\lambda\downarrow0$. This in turn gives that
\begin{equation}
B_\infty^\star(\lambda)=\sqrt\lambda
\bigl(c_\rho+o(1)\bigr),\qquad\lambda\downarrow0.
\end{equation}
Plugging this back in \eqref{Alimasymp} proves the desired claim.
\end{proofsect}

\subsection{Supercritical case}
\label{sec4.5}\noindent Here we will indicate the changes to the
arguments from the previous two sections that are needed to prove
Proposition~\ref{prop4.3}. We begin with an analogue of
Lemma~\ref{lemma4.6}:

\begin{lemma}
\label{lemma4.10} Let~$\rho,\rho'\in\MMsharp$ and define~$\rho_\alpha=(1-\alpha)\rho+\alpha\rho'$. Suppose that~$\zz(\rho)=\zzc$ and~$\zz(\rho_\alpha)>\zzc$ for all~$\alpha\in(0,1]$. Then for each~$\theta\ge1$, there is a constant~$K'(\theta)\in(0,\infty)$ such~that
\begin{equation}
\label{4.11b10} \limsup_{\alpha\downarrow0}
\frac{\BbbP_{\rho_\alpha}(|\eusmB^{(\theta)}|=\infty)}
{\zz(\rho_\alpha)-\zzc}\le K'(\theta).
\end{equation}
\end{lemma}

\begin{proofsect}{Proof}
The only important change compared to the proof of
Lemma~\ref{lemma4.6} is the derivation of the bound for~$\theta=1$. Indeed, in this case we use that~$\varkappa_{\rho_\alpha}(0)\ge
B_\infty^\star(0,\alpha)B_\infty(1,0,\alpha)$ in \eqref{LBbdsuper},
where~$B_\infty(1,0,\alpha)$ is the quantity~$B_\infty(\theta,\lambda)$
for~$\lambda=0$,~$\theta=1$ and~$\rho=\rho_\alpha$.
Applying $B_\infty^\star(0,\alpha)>0$ for all~$\alpha\in(0,1]$,
as follows by Theorem~\ref{thm3.1}(2), we find that
\eqref{4.11b10} holds with~$K'(1)=\frac{2b}{b-1}$. Once we set~$\lambda=0$, the rest of the proof can literally be copied.
\end{proofsect}

Next we need to state the appropriate version of
Lemma~\ref{lemma4.9a}:

\begin{lemma}
\label{lemma4.11} Let~$\rho,\rho'\in\MMsharp$ and define~$\rho_\alpha=(1-\alpha)\rho+\alpha\rho'$. Suppose that~$\zz(\rho)=\zzc$ and~$\zz(\rho_\alpha)>\zzc$ for all~$\alpha\in(0,1]$. Then
\begin{equation}
\label{4.11b11} \frac{\BbbP_{\rho_\alpha}(|\eusmB^{(\theta)}|=\infty)} {\zz(\rho_\alpha)-\zzc}=\psi_\rho(\theta)\,
\frac{\wE_\alpha(\BbbP_{\rho_\alpha}(|\eusmB^{(Q_\infty)}|=\infty))}
{\zz(\rho_\alpha)-\zzc}+o(1),\qquad \alpha\downarrow0,
\end{equation}
where~$\wE_\alpha$ is the
expectation corresponding to~$\wBbbP$ for measure~$\rho_\alpha$.
\end{lemma}

\begin{proofsect}{Proof}
Also in this case the required changes are only minuscule. First,
we have an analogue of~\eqref{Aasymp},
\begin{equation}
\label{4.11c11} \BbbP_{\rho_\alpha}\bigl(|\eusmB^{(\theta)}|=\infty\bigr)=b^nZ_n^{(\rho_\alpha)}(\theta)\,
\E_{n,\alpha}^{(\theta)}\bigl(\BbbP_{\rho_\alpha}\bigl(|\eusmB^{(Q_{n,1}^{(\theta)})}|=\infty\bigr)\bigr)
+\tilde\epsilon_n'(\alpha),
\end{equation}
where~$\E_{n,\alpha}^{(\theta)}$ is
the expectation~$\E_n^{(\theta)}$ and~$Z_n^{(\rho_\alpha)}$ the
object~$Z_n(\theta)$ for the underlying measure~$\rho_\alpha$ and
where~$\tilde\epsilon_n'(\alpha)$ is the quantity in \eqref{Aasymp}
for~$\lambda=0$ and~$\rho=\rho_\alpha$. We claim that
\begin{equation}
\label{limalpha}
\lim_{\alpha\downarrow0}\frac{\tilde\epsilon_n'(\alpha)}
{\zz(\rho_\alpha)-\zzc}=0
\end{equation}
for all finite~$n\ge1$.
Indeed, the entire derivation \twoeqref{UUn}{refbd} carries over, provided
we set~$\lambda=0$. The role of the ``small parameter'' is now taken over by~$\zz(\rho_\alpha)-\zzc$. A computation shows that~$\tilde\epsilon_n(\alpha)=O((\zz(\rho_\alpha)-\zzc)^2)$ as~$\alpha\downarrow0$, proving \eqref{limalpha}.

To finish the proof, it now remains to note that~$b^nZ_n^{(\rho_\alpha)}(\theta)\to b^n Z_n^{(\rho)}(\theta)$ as~$\alpha\downarrow0$ and that, by
Corollary~\ref{cor3.7} and the fact that~$\zz(\rho)=\zzc$, we have
$b^n Z_n^{(\rho)}(\theta)=\psi_\rho(\theta)+o(1)$ as~$n\to\infty$.
\end{proofsect}

Recall the definition of~$c_\rho$ in \eqref{crho}.
To prove Proposition~\ref{prop4.3}, we will need to know some basic continuity properties of~$c_\rho$ in~$\rho$. Note that these do not follow simply from the continuity of~$\alpha\mapsto\psi_{\rho_\alpha}(\theta)$, because also the expectation~$\wE$ in \eqref{crho} depends on the underlying measure.

\begin{lemma}
\label{lemma4.10a}
Let~$\rho,\rho'\in\MMsharp$ be such that~$\rho_\alpha=(1-\alpha)\rho+\alpha\rho'$ satisfies~$\zz(\rho_\alpha)>0$ for all~$\alpha\in[0,1]$. Let~$c_\rho$ be as in \eqref{crho}. Then~$\lim_{\alpha\downarrow0}c_{\rho_\alpha}=c_\rho$.
\end{lemma}

\begin{proofsect}{Proof}
Let~$\psi^*_{\rho_\alpha}(\theta)
=\E_{\rho_\alpha}(\psi_{\rho_\alpha}(X_\varnothing+\frac1b\theta))$. In general,~$\psi_{\rho_\alpha}(\theta)$ is Lipschitz continuous for~$\theta\ge1$. Thus,~$\psi_{\rho_\alpha}$ converges uniformly to~$\psi_\rho$ on compact sets of~$\theta$. Hence, we just need to show
\begin{equation}
\label{clmbd}
\lim_{\alpha\downarrow0}
\wE_\alpha\bigl(\psi^*_\rho(Q_\infty)^2\big|Q_\infty\ge1\bigr)
=\wE\bigl(\psi^*_\rho(Q_\infty)^2\big|Q_\infty\ge1\bigr).
\end{equation}
Choose~$n\ge1$ and replace~$\wE_\alpha$,~$\wE$ and~$Q_\infty$ by their finite-$n$ versions. By Corollary~\ref{Holder-cor}, the error thus incurred is uniformly small in~$\alpha\in[0,1]$. Hence, it is enough to show that
\begin{equation}
\label{finlim}
\lim_{\alpha\downarrow0}
\E_{n,\rho_\alpha}^{(\theta)}\bigl(\psi^*_\rho(Q_{n,1}^{(\theta)})^2
\big|Q_{n,1}^{(\theta)}\ge1\bigr)
=\E_{n,\rho}^{(\theta)}\bigl(\psi^*_\rho(Q_{n,1}^{(\theta)})^2
\big|Q_{n,1}^{(\theta)}\ge1\bigr),
\end{equation}
for some~$\theta\in[1,\theta_b]$, where~$\E_{n,\rho}^{(\theta)}$
denotes the expectation with respect to~$\BbbP_n^{(\theta)}$ for measure~$\rho$. However, in \eqref{finlim} only a finite number of coordinates are involved and the result follows.
\end{proofsect}

With Lemmas~\ref{lemma4.10},~\ref{lemma4.10a} and~\ref{lemma4.11}, we can finish the proof of Proposition~\ref{prop4.3}:

\begin{proofsect}{Proof of Proposition~\ref{prop4.3}}
From \eqref{4.11b11} we have
\begin{equation}
\frac{b-1}2\varkappa_{\rho_\alpha}(0)=B_\infty^\star(0,\alpha)^2
c_{\rho_\alpha}^{-2}+o\bigl(\zz(\rho_\alpha)-\zzc\bigr), \qquad\alpha\downarrow0.
\end{equation}
Using this in \eqref{LBbdsuper} and invoking 
Lemma~\ref{lemma4.10a}, we have
\begin{equation}
\frac{B_\infty^\star(0,\alpha)}
{\zz(\rho_\alpha)-\zzc}=bc_\rho^2+o(1),
\qquad\alpha\downarrow0.
\end{equation}
The proof is finished by plugging this back into \eqref{4.11b11} and invoking the
continuity of~$\alpha\mapsto\psi_{\rho_\alpha}(\theta)$.
\end{proofsect}

\section{Coupling argument}
\label{sec5}
\subsection{Coupling measure}
\label{sec5.1}\noindent The goal of this section is to define a
coupling of the measures~$\BbbP_n^{(\theta)}$ and~$\BbbP_n^{(\theta')}$ that appear in~\eqref{cbd}. As the first
step, we will write~$\BbbP_n^{(\theta)}(\cdot)$ as the
distribution of a time-inhomogeneous process. To have the process running in forward
time direction, we will need to express all quantities in terms of
the (more or less) original variables~$(X_k)$, which relate to the~$Y$'s through
\begin{equation}
\label{XYrel} 
X_k=Y_{n-k+1}\quad\text{or}\quad
Y_k=X_{n-k+1}, \qquad 1\le k\le n,
\end{equation}
see Section~\ref{sec3.3}. Abusing the notation slightly,~$\BbbP_n^{(\theta)}(\cdot)$ will temporarily be used to denote the
distribution of the~$X_1,\dots,X_n$ as well. We will return to the~$Y$'s in the proofs of Propositions~\ref{prop3.5}
and~\ref{prop3.6}.

Let~$Z_n(\theta)$ be as in \eqref{Zntheta} and note that, since~$\rho\in\MMsharp$, we have~$Z_n(\theta)>0$ for all~$n\ge0$ and all~$\theta\ge1$.
Given~$1\le k\le n-1$ and, for~$k>1$, a sequence~$(X_1,\dots, X_{k-1})\in[0,1]^{k-1}$,
let~$t_{n,k}^{(\theta)}(\cdot)=t_{n,k}^{(\theta)}(\,\cdot\,|X_1,\dots,X_{k-1})$
be given by
\begin{equation}
\label{tdef} t_{n,k}^{(\theta)}(x)=\frac{Z_{n-k-1}(x+\frac 1b
Q_{k-1}^{(\theta)})}{Z_{n-k}(Q_{k-1}^{(\theta)})}
\1_{\{Q_{k-1}^{(\theta)}\ge1\}},\qquad
0\le x\le1,
\end{equation}
where the indicator ensures that we are not dividing by zero. The~$(X_1,\dots, X_{k-1})$-dependence of~$t_{n,k}^{(\theta)}$ will be often left implicit. 

To interpret these objects, let us consider the case~$k=1$. Suppose 
that we wish to elucidate the distribution of~$X_1$ knowing that the process \textit{will} survive long enough to
produce an~$X_{n-1}$. (The variable~$X_n$ corresponds to~$Y_1$, which will be
uncorrelated with the other~$Y$'s.)
The only prior history we know is the value of~$\theta$; obviously
we are only interested in the case~$\theta\ge1$.
The total weight of all configurations is just~$Z_{n-1}(\theta)$; hence the denominator of \eqref{tdef}. Now, if~$X_1$ takes value~$x$, the weight of configurations in which the process survives is like the weight  of a string of length~$n-2$ with an effective ``$\theta$'' given by~$x+\frac1b\theta$. Hence~$Z_{n-2}(x+\frac1b\theta)$ in the numerator. (Notice that if~$x+\frac1b\theta<1$, this automatically vanishes.) We conclude that~$\BbbP_n^{(\theta)}(X_1\in\textd x)=t_{n,1}^{(\theta)}(x)\rho(\textd x)$. 

A similar reasoning shows that the probability of~$\{X_k\in\textd x\}$ given the values of~$X_1,\dots,X_{k-1}$ equals~$t_{n,k}^{(\theta)}(x)\rho(\textd x)$. This
allows us to view~$\BbbP_n^{(\theta)}$ as the distribution of an inhomogeneous process:

\begin{lemma}
\label{lemma5.1} For all~$\theta\ge1$, all~$n\ge1$ and all
Borel-measurable sets~$A\subset[0,1]^n$,
\begin{equation}
\BbbP_n^{(\theta)}(A)=\E\Bigl(\1_A \prod_{k=1}^{n-1}
t_{n,k}^{(\theta)}(X_k|X_1,\dots,X_{k-1})\Bigr).
\end{equation}
\end{lemma}

\begin{proofsect}{Proof}
The result immediately follows from the formula
\begin{equation}
\prod_{k=1}^{n-1} t_{n,k}^{(\theta)}(X_k|X_1,\dots,X_{k-1})=
\frac1{Z_{n-1}(\theta)}
\biggl\{\,\prod_{k=1}^{n-1}\1_{\{X_k+\frac 1b
Q_{k-1}^{(\theta)}\ge1\}}\biggr\},
\end{equation}
the identity~$Q_k^{(\theta)}=X_k+\frac 1b
Q_{k-1}^{(\theta)}$ and
the definition of~$\BbbP_n^{(\theta)}(\cdot)$, see
\eqref{Pntheta}.
\end{proofsect}

Next we will define the coupled measure. The idea is to use the
so-called Vasershtein coupling, see~\cite{Lindvall}, which
generates new (coupled) pairs from the ``maximal overlap'' of the
individual distributions. Let~$\theta,\theta'\ge1$, and suppose
that the corresponding sequences~$X=(X_1,\dots,X_{k-1})\in[0,1]^{k-1}$ and~$X'
=(X_1',\dots,X_{k-1}')\in[0,1]^{k-1}$ have been generated. Assume also
that a sequence~$(\omega_1,\dots,\omega_{k-1})\in\{0,1\}^{k-1}$ satisfying~$\omega_\ell\le\1_{\{X_\ell=X_\ell'\}}$ for all~$1\le\ell\le k-1$
has been generated. (This sequence marks down when~$X_\ell$ was
coupled with~$X_\ell'$. Note that we could have that~$X_\ell=X_\ell'$
even when~$X_\ell$ and~$X_\ell'$ are not coupled.) Let~$t$ be the
quantity~$t_{n,k}^{(\theta)}$ for the sequence~$X$ and let~$t'$ be the
corresponding quantity for the sequence~$X'$. Let
\begin{equation}
\label{RR}
R(\,\cdot\,) =R_{n,k}^{(\theta,\theta')}
(\,\cdot\,|X_1,\dots,X_{k-1};X_1',\dots,X_{k-1}';\omega_1,\dots,\omega_{k-1})
\end{equation}
be the transition kernel of the joint process, which is a
probability measure on~$[0,1]\times[0,1]\times\{0,1\}$ defined by
the expression
\begin{equation}
\label{doubleR} R\bigl(\textd x\!\times\!
\textd x'\!\times\!\{\omega\}\bigr)=\begin{cases} t(x)\wedge
t'(x)\,\rho(\textd  x)\delta_x(\textd x'),\quad&\text{if
}\omega=1,
\\
\frac1{1-q}\,{[t(x)-t'(x)]_+\,[t'(x')-t(x')]_+}\,\rho(\textd 
x)\rho(\textd  x'),\quad&\text{if }\omega=0.
\end{cases}
\end{equation}
Here~$t(x)\wedge t'(x)$ denotes the minimum of~$t(x)$ and~$t'(x)$
and~$[t(x)-t'(x)]_+$ denotes the positive part of~$t(x)-t'(x)$.
The quantity~$q=q_{n,k;X,X'}^{(\theta,\theta')}$ is given by
\begin{equation}
\label{varpidef} 
q=\int t(x)\wedge t'(x)\,\rho(\textd x)=
1-\int \bigl[t(x)-t'(x)\bigr]_+\rho(\textd x).
\end{equation}
The interpretation of \eqref{doubleR} is simple: In order to
sample a new triple~$(X_k,X'_k,\omega_k)$, 
we first choose~$\omega_k\in\{0,1\}$ with~$\Prob(\omega_k=1)=q$. If~$\omega_k=1$, the pair~$(X_k,X'_k)$ is sampled from distribution~$\frac1q
t(x)\wedge t'(x)\,\rho(\textd  x)\delta_x(\textd x')$---and, in
particular,~$X_k$ gets glued together with~$X_k'$---while
for the case~$\omega_k=0$ we use the distribution in the
second line of \eqref{doubleR}.

\begin{remark}
It turns out that whenever the above processes~$X$ and~$X'$ 
have glued together, they
have a tendency to stay glued. However, the above coupling is
\textit{not} monotone, because the processes may come apart
no matter how long they have been glued together. Our strategy lies
in showing that~$q$ tends to one rapidly enough so that the number
of ``unglueing'' instances is finite almost surely.
\end{remark}

Let~$\BbbP_n^{(\theta,\theta')}(\cdot)$ be the probability measure
on~$[0,1]^n\times[0,1]^n\times\{0,1\}^n$ assigning mass
\begin{equation}
\label{doubleP}
\BbbP_n^{(\theta,\theta')}(B)=\sum_{(\omega_k)}\int_B
\rho(\textd x_n)\rho(\textd x_n')\1_{\{\omega_n=1\}}
\prod_{k=1}^{n-1}
R_{n,k;x,x',\omega}^{(\theta,\theta')}\bigl(\textd x_k\!\times\!
\textd x_k'\!\times\!\{\omega_k\}\bigr)
\end{equation}
to any Borel-measurable set~$B\subset[0,1]^n\times[0,1]^n\times\{0,1\}^n$. Here~$R_{n,k;x,x',\omega}^{(\theta,\theta')} (\textd x_k\!\times\!
\textd x_k'\!\times\!\{\omega_k\})=R_{n,k}^{(\theta,\theta')}
(\textd x_k\!\times\! \textd x_k'\!\times\!\{\omega_k\}
|x_1,\dots,x_{k-1};x_1',\dots,x_{k-1}';\omega_1,\dots,\omega_{k-1})$.
As can be expected from the construction,~$\BbbP_n^{(\theta)}(\cdot)$ and~$\BbbP_n^{(\theta')}(\cdot)$ are
the first and second marginals of~$\BbbP_n^{(\theta,\theta')}(\cdot)$, respectively:

\vbox{
\begin{lemma}
\label{lemma5.2} Let~$\theta,\theta'\ge1$. Then
\begin{align}
\label{marginal1} \BbbP_n^{(\theta,\theta')}\bigl(A\times
[0,1]^n\times\{0,1\}^n\bigr)&=\BbbP_n^{(\theta)}(A)\\
\intertext{and} \label{marginal2}
\BbbP_n^{(\theta,\theta')}\bigl(\,[0,1]^n\times
A\times\{0,1\}^n\bigr) &=\BbbP_n^{(\theta')}(A),
\end{align}
for all Borel-measurable~$A\subset[0,1]^n$.
\end{lemma}
}

\begin{proofsect}{Proof}
To prove formula \eqref{marginal1}, let~$X=(X_1,\dots,X_{k-1})$ and~$X'=(X_1',\dots,X_{k-1}')$ be two sequences from~$[0,1]^{k-1}$. If~$Q_{k-1}^{(\theta)}\ge1$ and the same holds for the corresponding
quantity for the sequence~$X'$, let~$t(\cdot)=t_{n,k}^{(\theta)}(\cdot)$,~$t'(\cdot)=t_{n,k}^{(\theta')}(\cdot)$, and let~$R(\cdot)$ and~$q$
be as in \eqref{doubleR} and \eqref{varpidef}, respectively. Using 
\eqref{varpidef} we have, for all Borel sets~$C\subset[0,1]$,
\begin{equation}
\label{cpl} \sum_{\omega\in\{0,1\}}\int_{C\times[0,1]}\,
\!\!\!\!\!R(\textd x\!\times\!\textd x'\!\times\!\{\omega\})=
\int_C \bigl(t(x)\wedge t'(x)+[t(x)-t'(x)]_+\bigr)\rho(\textd 
x)= \int_C t(x)\rho(\textd x).
\end{equation}
In other words, 
the first marginal of the coupled process is a process on~$[0,1]$ with the transition kernel~$t(\cdot)\rho(\cdot)$, which,
as shown in Lemma~\ref{lemma5.1}, generates~$\BbbP_n^{(\theta)}$.
This proves \eqref{marginal1}; the proof of \eqref{marginal2} is
analogous.
\end{proofsect}

Clearly, the number~$q$ represents the probability that the two
processes get coupled. The following lemma provides a bound that
will be useful in controlling~$q$:

\begin{lemma}
\label{lemma5.3a} Let~$\theta,\theta'\ge1$,~$1\le k\le n-1$ and~$X=(X_1,\dots,X_{k-1})\in[0,1]^{k-1}$ and~$X'=(X_1',\dots,X_{k-1}')\in[0,1]^{k-1}$. Let~$Q$ be the quantity~$Q_{k-1}^{(\theta)}$ corresponding to~$X$ and let~$Q'$ be the quantity~$Q_{k-1}^{(\theta')}$ 
corresponding to~$X'$. If~$Q\wedge Q'\ge1$, then
\begin{equation}
\label{varpibd}
q_{n,k;X,X'}^{(\theta,\theta')}\ge\frac{Z_{n-k}(Q\wedge
Q')}{Z_{n-k}(Q\vee Q')}.
\end{equation}
\end{lemma}

\begin{proofsect}{Proof}
Let~$t$ be the
quantity~$t_{n,k}^{(\theta)}$ for the sequence~$X$ and let~$t'$ be the
corresponding quantity for the sequence~$X'$. By
inspection of \eqref{tdef} and monotonicity of~$\theta\mapsto
Z_n(\theta)$,
\begin{equation}
t(x)\ge\frac{Z_{n-k-1}(x+\frac1b(Q\wedge Q'))} {Z_{n-k}(Q\vee
Q')},
\end{equation}
and similarly for~$t'(x)$. From here the claim follows by
integrating with respect to~$\rho(\textd x)$.
\end{proofsect}

\subsection{Domination by a discrete process}
\label{sec5.2}\noindent The goal of this section is to show that
the coupled measure defined in the previous section has the
desirable property that, after a finite number of steps, the
processes~$X$ and~$X'$ get stuck forever. Since the
information about coalescence of~$X$ and~$X'$ is encoded into the
sequence~$\omega$, we just need to show that, eventually,~$\omega_k=1$. For technical reasons, we will concentrate from the
start on infinite sequences~$(\omega_k)_{k\in\N}$: Let~$P_n^{(\theta,\theta')}(\cdot)$ be the law of~$(\omega_k)_{k\in\N}\in\{0,1\}^\N$ induced by the distribution~$\BbbP_n^{(\theta,\theta')}(\cdot)$ and the requirement~$P_n^{(\theta,\theta')}(\omega_k=1,\,k\ge n)=1$.

The coalescence of~$X$ and~$X'$ will be shown by a comparison with
a simpler stochastic process on~$\{0,1\}^\N$ whose law will be
distributionally lower than~$P_n^{(\theta,\theta')}(\cdot)$, i.e., in the
FKG~sense. Let~$\preccurlyeq$ be the partial order on~$\omega,\omega'\in\{0,1\}^\N$ defined by
\begin{equation}
\omega\preccurlyeq\omega'\quad\Leftrightarrow\quad
\omega_k\le\omega_k',\qquad k\ge1.
\end{equation}
Next, note that, by~$x_\star>\frac{b-1}b$, we have~$1-b(1-x_\star)>\theta_b-1$. Choose a number~$\delta_\rho\in(\theta_b-1,1-b(1-x_\star))$ and, noting 
that~$\rho([1-\frac {1-\delta_\rho}
b,x_\star])>0$, define a collection of weights~$(\lambda_\rho(s))$ by 
\begin{equation}
\label{lambdan} \frac{1-\lambda_\rho(s)}{\lambda_\rho(s)}
=\sum_{k\ge s}\,\,\sup_{\theta-\theta'\le\delta_\rho b^{-k}}
\frac {\rho\bigl([1-\frac \theta b,1-\frac{\theta'}b)\bigr)}
{\rho\bigl([1-\frac {1-\delta_\rho} b,x_\star]\bigr)}, \qquad
s\in\N\cup\{0\}.
\end{equation}
Note that~$s\mapsto\lambda(s)$ is increasing. It is also easy to
verify that~$\lambda_\rho(\cdot)\in(0,1]$, 
so any of these weights can
be interpreted as a probability. 
This allows us to define a process on~$(\omega'_k)_{k\in\N}\in\{0,1\}^\N$, with the transition kernel
\begin{equation}
\label{lambdap}
p_\rho(\,\omega'_k=1\,|\,\omega'_1,\dots,\omega'_{k-1})
=\lambda_\rho\bigl(\min\{0\le j\le k-1\colon
\omega_{k-j-1}'=0\}\bigr),
\end{equation}
where, for definiteness,
we set~$\omega_0'=0$. Let~$\widetilde P_\rho(\cdot)$ denote
the law of the entire process with transition probabilities~$p_\rho(\,\cdot\,|\,\cdot\,)$ and ``initial'' value~$\omega_0'=0$.

\begin{proposition}
\label{prop5.3} Let~$\rho\in\MMsharp$ and let~$\delta_\rho$ be as
above. For all~$n\ge1$ and all~$\theta,\theta'$ with~$1\le\theta,\theta'\le\theta_b$, the
measure~$P_n^{(\theta,\theta')}(\cdot)$ stochastically
dominates~$\widetilde P_\rho(\cdot)$ in partial
order~$\preccurlyeq$.
\end{proposition}

Let~$\delta_\rho$ be fixed for the rest of this Subsection. In order to
give a proof of Proposition~\ref{prop5.3}, we first
establish a few simple bounds.

\begin{lemma}
\label{lemma5.3b} Let~$\rho\in\MMsharp$ and let~$\delta_\rho$ be as above. 
Let~$n\ge0$ and suppose~$\theta,\theta'\ge1$ satisfy~$0\le\theta-\theta'\le\delta_\rho
b^{-k}$ for some~$k\ge0$. Then
\begin{equation}
\frac{Z_n(\theta')}{Z_n(\theta)}\ge\lambda_\rho(k).
\end{equation}
\end{lemma}

\begin{proofsect}{Proof}
Consider a configuration~$X_1,\dots,X_n$ which contributes to~$Z_n(\theta)$ but \textit{not} to~$Z_n(\theta')$. This implies that there is an~$\ell\in\{1,\dots,n\}$ where~$Q_\ell^{(\theta)}\ge1$ but~$Q_\ell^{(\theta')}<1$. With this in mind, we
claim the identity
\begin{equation}
\prod_{m=1}^n\1_{\{Q_m^{(\theta)}\ge1\}}
-\prod_{m=1}^n\1_{\{Q_m^{(\theta')}\ge1\}}
=\sum_{\ell=1}^n
\biggl[\,\prod_{m=1}^{\ell-1}\1_{\{Q_m^{(\theta')}\ge1\}}\biggl]
\,\1_{\{Q_\ell^{(\theta')}<1\le
Q_\ell^{(\theta)}\}}\,
\biggl[\,\prod_{m=\ell+1}^n\1_{\{Q_m^{(\theta)}\ge1\}}\biggr].
\end{equation}
Thence,
\begin{equation}
Z_n(\theta)-Z_n(\theta')=\sum_{\ell=1}^n\E\biggl(
Z_{n-\ell}\bigl(Q_\ell^{(\theta)}\bigr)\,\1_{\{Q_\ell^{(\theta')}<1\le
Q_\ell^{(\theta)}\}}\,\prod_{m=1}^{\ell-1}\1_{\{Q_m^{(\theta')}\ge1\}}\biggr).
\end{equation}
Since~$\theta-\theta'\le\delta_\rho b^{-k}$, we have~$Q_\ell^{(\theta)}-1\le Q_\ell^{(\theta)}-Q_\ell^{(\theta')}\le
\delta_\rho b^{-k-\ell}$ for any~$\ell$ contributing on the
right-hand side. In particular, we have~$Q_\ell^{(\theta)}\le
1+\tfrac{\delta_\rho}b$, which implies~$Z_{n-\ell}(Q_\ell^{(\theta)})\le Z_{n-\ell}(1+\tfrac{\delta_\rho}b)$. Then
\begin{equation}
\label{one} Z_n(\theta)-Z_n(\theta')\le\sum_{\ell=1}^n
Z_{n-\ell}\bigl(1+\tfrac{\delta_\rho}b\bigr) \E\biggl(
\rho\bigl([1-\tfrac1bQ_{\ell-1}^{(\theta)},1-\tfrac1bQ_{\ell-1}^{(\theta')})\bigr)
\,
\prod_{m=1}^{\ell-1}\1_{\{Q_m^{(\theta')}\ge1\}}\biggr),
\end{equation}
or, replacing~$\rho([1-\tfrac1bQ_{\ell-1}^{(\theta)},1-\tfrac1bQ_{\ell-1}^{(\theta')}))$ by its maximal value,
\begin{equation}
\label{one1} Z_n(\theta)-Z_n(\theta')\le\sum_{\ell=1}^n
Z_{n-\ell}\bigl(1+\tfrac{\delta_\rho}b\bigr)\,
Z_{\ell-1}(\theta')
\sup_{\vartheta-\vartheta'\le\delta_\rho b^{-k-\ell+1}}
\!\!\!\!\rho\bigl([1-\tfrac \vartheta b,1-\tfrac{\vartheta'}b)\bigr).
\end{equation}
On the other hand, by simply demanding that~$X_\ell\ge1-\tfrac{1-\delta_\rho}b$ (which implies~$Q_\ell^{(\theta)}\ge1+\tfrac{\delta_\rho}b$) in \eqref{ZZrel} we
have for all~$1\le\ell\le n$ that
\begin{equation}
\label{two} Z_n(\theta')\ge
Z_{n-\ell}\bigl(1+\tfrac{\delta_\rho}b\bigr)
\rho\bigl([1-\tfrac{1-\delta_\rho}b,x_\star]\bigr)Z_{\ell-1}(\theta').
\end{equation}
Using \eqref{two} in \eqref{one1}, and applying \eqref{lambdan},
we have
\begin{equation}
\label{three}
Z_n(\theta)-Z_n(\theta')\le
\frac{1-\lambda_\rho(k)}{\lambda_\rho(k)}Z_n(\theta'),
\end{equation}
whereby the claim directly follows.
\end{proofsect}

Next we prove a bound between kernels \eqref{doubleR} and
\eqref{lambdap}:

\begin{lemma}
\label{lemma5.4} Let~$1\le k\le n-1$ and let~$\omega'=(\omega_1',\dots,\omega_{k-1}')\in\{0,1\}^{k-1}$,~$X=(X_1,\dots,X_{k-1})\in[0,1]^{k-1}$,~$X'= (X_1',\dots,X_{k-1}')\in[0,1]^{k-1}$
and~$\omega=(\omega_1,\dots,\omega_{k-1})\in\{0,1\}^{k-1}$. 
For all~$\theta,\theta'\ge1$ and all~$\ell=1,\dots, k-1$, let
$Q_\ell^{(\theta)}$ correspond to~$X$ via \eqref{Qntheta}, 
and let~$Q_\ell^{(\theta')}$ correspond to~$X'$. Suppose that
\begin{equation}
\label{couplcond}
Q_j^{(\theta)}\ge1,\quad Q_j^{(\theta')}\ge1\quad
\text{and}\quad\omega_j'\le\omega_j \le\1_{\{X_j=X_j'\}},\qquad
j=1,\dots,k-1.
\end{equation}
If~$R_{n,k;X,X',\omega}^{(\theta,\theta')}(\cdot)$
is the quantity defined in~\eqref{doubleP}, then
\begin{equation}
\label{ineq} 
R_{n,k;X,X',\omega}^{(\theta,\theta')}
\bigl(\{\omega_k=1\}\bigr) \ge
p_\rho(\,\omega_k'=1\,|\,\omega_1',\dots,\omega_{k-1}'),
\end{equation}
for all~$\theta,\theta'$ with~$1\le\theta,\theta'\le\theta_b$.
\end{lemma}

\begin{proofsect}{Proof}
Note that, since~$1\le\theta,\theta'\le\theta_b$ and~$1+\delta_\rho\ge\theta_b$, we have~$1\le Q_\ell^{(\theta)},Q_\ell^{(\theta')}\le 1+\delta_\rho$ and thus~$|Q_\ell^{(\theta)}-Q_\ell^{(\theta')}|\le\delta_\rho$ for all~$\ell=1,\dots,k-1$. This allows us to define the quantity
\begin{equation}
s=\max\bigl\{\ell\colon 0\le\ell\le k,\,|Q_{k-1}^{(\theta)}-Q_{k-1}^{(\theta')}|\le\delta_\rho b^{-\ell}\bigr\}.
\end{equation} 
By Lemmas~\ref{lemma5.3a}
and~\ref{lemma5.3b}, we have~$R(\{\omega_k=1\}) \ge\lambda_\rho(s)$, where~$R(\cdot)$ stands for the quantity on the left-hand side of \eqref{ineq}.
Recall our convention~$\omega_0'=0$ and let
\begin{equation}
\label{5.27}
s'=\min\bigl\{0\le j\le k-1\colon\omega_{k-j-1}'=0\bigr\}.
\end{equation}
In other words,~$s'$ is the length of the largest contingent block of~$1$'s in~$\omega'$
directly preceding~$\omega_k'$.
We claim that~$s\ge s'$. Indeed, by our previous reasoning,
$|Q_{k-s'-1}^{(\theta)}-Q_{k-s'-1}^{(\theta')}|\le\delta_\rho$. By our assumptions,
$1=\omega_j'\le\1_{\{X_j=X_j'\}}$ and, therefore,~$X_j=X_j'$ 
for all~$j=k-s',\dots,k-1$. This~implies
\begin{equation}
\bigl|Q_{k-1}^{(\theta)}-Q_{k-1}^{(\theta')}\bigr|\le \delta_\rho b^{-s'}
\end{equation}
and hence~$s\ge s'$. Using that~$s'$ is the argument of~$\lambda$ in \eqref{lambdap} we have
$R(\{\omega_k=1\})\ge\lambda_\rho(s)\ge \lambda_\rho(s')=
p_\rho(\,\omega_k'=1\,|\,\omega_1',\dots,\omega_{k-1}')$. This proves the claim.
\end{proofsect}

Now we are ready to prove Proposition~\ref{prop5.3}:

\begin{proofsect}{Proof of Proposition~\ref{prop5.3}}
The inequality \eqref{ineq} is a sufficient
condition for the existence of so-called Strassen's coupling,
see~\cite{Lindvall}. In particular, the inhomogeneous-time 
process generating the triples~$(X_k,X_k',\omega_k)$ can be coupled with
the process generating~$\omega_k'$ in such a way that \eqref{couplcond} holds at
all times less than~$n$. The~$(\omega,\omega')$ marginal of this process will be, by definition, concentrated on~$\{\omega\succcurlyeq\omega'\}$. Since~$\omega_k=1$ for~$k>n$,~$P_n^{(\theta,\theta')}$-almost surely, the required stochastic domination follows.
\end{proofsect}

\subsection{Existence of the limiting measure}
\label{sec5.3}\noindent The goal of this section is to show that,
under proper conditions, the process~$\omega'$ with distribution~$\widetilde P_\rho(\cdot)$ equals one except at a finite number of
sites. Then we will give the proof of Proposition~\ref{prop3.5}.~Let
\begin{equation}
\label{pn}
p_n=\begin{cases}\bigl(1-\lambda_\rho(n)\bigr)\prod_{k=0}^{n-
1}\lambda_\rho(k),\qquad&
\text{if }n\in\N\cup\{0\},\\
\prod_{k=0}^\infty\lambda_\rho(k),\qquad& \text{if }n=\infty,
\end{cases}
\end{equation}
and observe that~$p_n$ is the probability of seeing a block of~$1$'s
of length~$n$ in the prime configuration.
We begin with an estimate of~$\lambda(k)$:

\begin{lemma}
\label{lemma5.8a} For each~$\rho\in\MMsharp$, there is~$C(\rho)<\infty$ and~$\varpi>0$ such that
\begin{equation}
1-\lambda_\rho(k)\le C(\rho)e^{-\varpi k}.
\end{equation}
Moreover, the quantity~$C(\rho)$ is bounded away from infinity
uniformly in any subset~$\NN\subset\MMsharp$ with finitely many
extreme points.
\end{lemma}

\begin{proofsect}{Proof}
Let~$\phi_\rho$ be the density of~$\rho$ with respect to the Lebesgue measure on~$[0,1]$. Then
\begin{equation}
\label{5.31} 
\sup_{\theta-\theta'\le\delta_\rho b^{-n}}\,
\rho\bigl([1-\tfrac \theta b,1-\tfrac{\theta'}b)\bigr)\le
\delta_\rho b^{-n}\Vert\phi_\rho\Vert_\infty.
\end{equation}
The claim then follows by inspection of \eqref{lambdan} with~$\varpi=\log b$ and an appropriate choice
of~$C(\rho)$. The bound on~$C(\rho)$ is
uniform in any~$\NN$ with the above properties, because the bound~$\Vert\phi_\rho\Vert_p<\infty$ is itself uniform.
\end{proofsect}

The preceding estimate demonstrates that the discrete process locks, and
in fact does so fairly rapidly. Indeed, we now have~$p_\infty>0$, which
ensures that eventually the configuration is all ones, and further that the
$p_n$ tend to zero exponentially. It remains to show that the  waiting times
till locking are themselves exponential.

\begin{lemma}
\label{lemma5.5} Let~$\rho\in\MMsharp$ and, for~$n\ge1$, let~$\EE(n)=\{\omega'\in\{0,1\}^\N\colon \omega_j'=1,\, j\ge n\}$.
Let~$\alpha_0>0$ be such that~$\varphi(\alpha)=\sum_{0\le k<\infty}e^{\alpha (k+1)}p_k <\infty$
for all~$\alpha\in(0,\alpha_0)$. Then
\begin{equation}
\label{bdonEE} \widetilde P_\rho\bigl(\EE(n)^{\text{\rm c}}\bigr)
\le n\, e^{-\mu(\rho) n}, \qquad n\ge1,
\end{equation}
where
\begin{equation}
\mu(\rho)=\sup
\bigl\{\alpha\ge0\colon\varphi(\alpha)\le1\bigr\}.
\end{equation}
\end{lemma}

We note that both quantities~$\alpha_0$ and~$\mu(\rho)$ are nontrivial. Indeed,~$\alpha_0\ge\varpi>0$ and, since~$p_\infty$ can be written as~$p_\infty=1-\sum_{n\ge0}p_n>0$, we have that~$\mu(\rho)>0$.

\begin{proofsect}{Proof}
An inspection of \eqref{lambdap} shows that ``blocks of 1's'' form
a renewal process. Indeed, suppose~$\xi_\ell$ for~$\ell=1,\dots,k-1$ mark down the lengths of first~$k-1$ ``blocks
of 1's'' including the terminating zero (i.e.,~$\xi_\ell=n$ refers
to a block of~$n-1$ ones and followed by a zero). Denoting~$N_{k-1}=\sum_{j=1}^{k-1}\xi_j$, the~$k$-th block's length is then
\begin{equation}
\xi_k=\min\{j>0\colon \omega_{j+N_{k-1}}'=0\}.
\end{equation}
As is seen from \eqref{lambdap},~$(\xi_\ell)$ can be continued
into an infinite sequence of i.i.d.\ random variables on~$\N\cup\{\infty\}$ with distribution~$\Prob(\xi_k=n+1)=p_n$, where~$p_n$ is as in \eqref{pn}. The physical sequence terminates after
the first~$\xi_k=\infty$ is encountered. Let~$\GG_n(k)$ be the
event that~$\xi_1,\dots,\xi_k$ are all finite and~$\sum_{i=1}^k\xi_i> n$. Then, clearly,~$\EE(n)^{\text{\rm
c}}=\bigcup_{k=1}^n\GG_n(k)$.

The probability of~$\GG_n(k)$ is easily bounded using the
exponential Chebyshev inequality:
\begin{equation}
\Prob\bigl(\GG_n(k)\bigr)\le \varphi(\alpha)^ke^{-\alpha n},
\qquad 0\le\alpha<\alpha_0.
\end{equation}
Noting that~$\sum_{k=1}^n\varphi(\alpha)^k\le n$ for~$\alpha\le\mu(\rho)$, the claim follows.
\end{proofsect}

Now we are finally ready to prove Proposition~\ref{prop3.5}:

\begin{proofsect}{Proof of Proposition~\ref{prop3.5}}
Let~$\rho\in\MMsharp$ and~$n$ be fixed. Let~$k\le n$ and suppose 
that~$f$ is a function that depends only on the first~$k$ of the~$Y$-coordinates. 
Let~$\theta_0>\theta_b$ and let~$\theta,\theta'\in[1,\theta_0]$. 
Noting
that~$\BbbP_n^{(\theta)}(\cdot|Q_{n,m}^{(\theta)}\in
\textd Q)=\BbbP_{n-m}^{(Q)}(\cdot)$, we have
\begin{equation}
\label{odhad}
\bigl|\E_{n+1}^{(\theta)}(f)-\E_n^{(\theta')}(f)\bigr|
\le\,\E_{n+1}^{(\theta)}
\Bigl(\bigl|\E_{n}^{(Q_{n+1,n}^{(\theta)})}(f)
-\E_{n}^{(\theta')}(f)\bigr|\Bigr).
\end{equation} 
Since~$Q_{n+1,n}^{(\theta)}\in[1,\theta_0]$ by our choice of~$\theta$,
we just need to estimate~$|\E_n^{(\theta)}(f)
-\E_n^{(\theta')}(f)|$ by the right-hand side of \eqref{cbd} for
all~$\theta,\theta'\in[1,\theta_0]$. 

Introduce the quantity
\begin{equation}
D_n(f)=\sup\bigl\{|\E_n^{(\theta)}(f)-\E_n^{(\theta')}(f)|\colon
\theta,\theta'\in[1,\theta_0]\bigr\}.
\end{equation}
We need to show~$D_n(f)$ is exponentially small in~$n$.
By Lemmas~\ref{lemma5.1}, \ref{lemma5.2}, and
Proposition~\ref{prop5.3}, the probability that~$X_i\ne X_i'$ for
some~$n-k\le i\le n$ under the coupling measure~$\BbbP_n^{(\theta,\theta')}(\cdot)$ is dominated by the
probability that~$\omega_i'=0$ for some~$n-k\le i\le n$ under~$\widetilde P_\rho(\cdot)$. Since~$f$ depends only on the first~$k$ 
of the~$Y$ variables (i.e., the \textit{last}~$k$ of the~$X$ variables), 
the coupling inequality gives us
\begin{equation}
\label{couplineq}
\bigl|\E_n^{(\theta)}(f)-\E_n^{(\theta')}(f)\bigr| \le2\Vert
f\Vert_\infty\, \widetilde P_\rho\bigl(\EE(n-k)^{\text{\rm
c}}\bigr),
\end{equation}
where~$\EE(n-k)$ is as in Lemma~\ref{lemma5.5}.

Let~$\mu=\mu(\rho)$ be as in Lemma~\ref{lemma5.5}. Then
\eqref{bdonEE} and \eqref{couplineq} give
\begin{equation}
D_n(f)\le 2\Vert f\Vert_\infty(n-k)\, e^{-\mu(n-k)}\le 4(\mu
e)^{-1} \Vert f\Vert_\infty e^{-\frac12\mu(n-k)},
\end{equation}
This proves \eqref{cbd} with~$\zeta=\frac12\mu$ 
and~$A=4(\mu e)^{-1}$. The bounds~$\zeta>0$ and~$A<\infty$ are uniform in sets~$\NN\subset\MMsharp$
with finitely-many extreme points, because the bound~$\mu(\rho)>0$ is itself uniform. 
The existence of the limit~\eqref{limeq} and its independence of~$\theta$ 
is then a direct consequence of~\eqref{cbd}.
\end{proofsect}

\subsection{Distributional identity}
\label{sec5.4}\noindent Here we will show the validity of the
distributional identity \eqref{disteq}. The proof we follow
requires establishing that the distribution of~$Q_\infty$ has no
atom at~$Q_\infty=1$:

\begin{lemma}
\label{lemmacont} Let~$\rho\in\MMsharp$. Then~$\wBbbP(Q_\infty=1)=0$.
\end{lemma}

\begin{proofsect}{Proof}
Notice that the almost-sure bound~$Q_{n,1}^{(1)}\le Q_\infty\le
Q_{n,1}^{(\theta_b)}$ holds for all~$n\ge1$, with~$Q_{n,1}^{(1)}\uparrow Q_\infty$ and~$Q_{n,1}^{(\theta_b)}\downarrow Q_\infty$ as~$n\to\infty$.
Therefore,
\begin{equation}
\wBbbP(Q_\infty=1)=\lim_{n\to\infty}
\wBbbP\bigl(Q_{n,1}^{(1)}<1,\, Q_{n,1}^{(\theta_b)}\ge1\bigr).
\end{equation}
But~$Y_1$ is unconstrained under~$\wBbbP(\cdot)$ which by~$0\le
Q_{n,1}^{(\theta_b)}-Q_{n,1}^{(1)}\le(\theta_b-1)b^{-n}$ allows us
to write
\begin{equation}
\label{QQdiff} \wBbbP\bigl(Q_{n,1}^{(1)}<1,\,
Q_{n,1}^{(\theta_b)}\ge1\bigr)\le\text{ l.h.s. of~\eqref{5.31}}.
\end{equation}
Hence,~$\wBbbP(Q_{n,1}^{(1)}<1,\, Q_{n,1}^{(\theta_b)}\ge1)\to0$
as~$n\to\infty$ and we have~$\wBbbP(Q_\infty=1)=0$, as claimed.
\end{proofsect}

\begin{proofsect}{Proof of Proposition~\ref{prop3.6}}
Let~$X$ be a random variable with distribution~$\BbbP(\cdot)=\rho(\cdot)$, independent of~$Y_1,Y_2,\dots$, and
let~$\theta\ge1$. For all~$a\in\R$, define the (distribution)
functions
\begin{equation}
F_n^{(\theta)}(a)=\BbbP_n^{(\theta)}\bigl(Q_{n,1}^{(\theta)}\ge
a\bigr).
\end{equation}
and
\begin{equation}
\widetilde F_n^{(\theta)}(a)= \BbbP\otimes\BbbP_n^{(\theta)}
\biggl(X+\frac{Q_{n,1}^{(\theta)}}b\ge a
,\,Q_{n,1}^{(\theta)}\ge1\biggr).
\end{equation}
Since~$Q_{n,1}^{(\theta)}\overset\DD=Q_{n+1,2}^{(\theta)}$,~$X\overset\DD=Y_1$ and~$Y_1+\frac1bQ_{n+1,2}^{(\theta)}=Q_{n+1,1}^{(\theta)}$, these
functions obey the relation
\begin{equation}
\label{finrel} \widetilde F_n^{(\theta)}(a)
=F_n^{(\theta)}(1)\,F_{n+1}^{(\theta)}(a), \qquad n\ge1,\,a\in\R.
\end{equation}
Let~$F(a)=\wBbbP(Q_\infty\ge a)$ and let
\begin{equation}
\widetilde F(a)=\BbbP\otimes\wBbbP \biggl(X+\frac{Q_\infty}b\ge a
,\,Q_\infty\ge1\biggr).
\end{equation}
Both~$F(\cdot)$ and~$\widetilde F(\cdot)$ are non-increasing,
left-continuous and they both have a right-limit at every~$a\in\R$. In particular, both functions are determined by their
restriction to any dense subset of~$\R$. The proof then boils down
to showing that there is a set~$A\subset\R$ dense in~$\R$ such
that
\begin{alignat}{2}
\label{toprove}
\lim_{n\to\infty} F_n^{(\theta)}(a)&=F(a)\qquad &a\in A\cup\{1\},\\
\intertext{and} \label{alsotoprove}\lim_{n\to\infty}\widetilde
F_n^{(\theta)}(a)&=\widetilde F(a),\qquad &a\in A.
\end{alignat}
Indeed, then \eqref{finrel} implies~$\widetilde F(a)=F(1)F(a)$ for
all~$a\in A$, which by continuity extends to all~$a\in\R$, proving
\eqref{disteq}. 

Let~$A$ be the set of continuity points of both~$F(\cdot)$ and~$\widetilde F(\cdot)$. Clearly,~$A^{\text{c}}$ is countable and
hence~$A$ is dense in~$\R$. The limits in \eqref{toprove} will be
taken in too stages; first we take the limit of the distribution
and then that of the event. Since $Q_{m,1}^{(1)}\le Q_{n,1}^{(\theta)}\le
Q_{m,1}^{(\theta_b)}$ for any~$m\le n$, we have, by \eqref{limeq},
\begin{equation}
\wBbbP (Q_{m,1}^{(1)}\ge a )\le \liminf_{n\to\infty}F_n^{(\theta)}(a)
\le\limsup_{n\to\infty}F_n^{(\theta)}(a)\le
\wBbbP(Q_{m,1}^{(\theta_b)}\ge a)
\end{equation}
for all~$\theta\ge1$ and all~$m\ge1$. The~$m\to\infty$ of the extremes
exists by monotonicity. Since $Q_{m,1}^{(\theta_b)}\ge Q_\infty$,
the right-hand side converges to~$F(a)$. As for the
left-hand side, it is clear that the event $\{Q_\infty>a\}$ implies that,
\textit{eventually}, $\{Q_{m,1}^{(1)}\ge a\}$ occurs. Thus the limit of the extreme left is at least as big as~$\wBbbP(Q_\infty>a)$.
However, the latter equals~$F(a)$ because, by assumption,~$a$ is a
continuity point of~$F$. 
This proves \eqref{toprove}. The argument
for the limit \eqref{alsotoprove} is fairly similar; the
right-hand side will directly converge to~$\widetilde F(a)$, while
the limit of the left hand side will be no smaller than
$\BbbP\otimes\wBbbP(X+\frac1bQ_\infty>a,\,Q_\infty>1)$. However, by
Lemma~\ref{lemmacont} we have that~$\wBbbP(Q_\infty=1)=0$ and thus
the limit equals~$\widetilde F(a)$, because~$a\in A$.
\end{proofsect}

\begin{proofsect}{Proof of Corollary~\ref{Holder-cor}}
The proof of~$Q_{n,1}^{(\theta)}\overset\DD\longrightarrow Q_\infty$ 
is immediate from
\eqref{toprove}. To prove \eqref{Hbd}, we note that \eqref{HBBbd}
and \eqref{Qntheta} imply the
deterministic bounds
\begin{equation}
\label{eins}
\bigl|f(Q_{2n,1}^{(\theta)})-f(Q_{n,1}^{(\theta_b)})\bigr|\le C\Vert f\Vert_\infty\,
b^{-n}\theta_0,
\end{equation}
and
\begin{equation}
\label{zwei}
\bigl|f(Q_\infty)-f(Q_{n,1}^{(\theta_b)})\bigr|\le C\Vert f\Vert_\infty\,
b^{-n}\theta_0,
\end{equation}
where we used that~$Q_{2n,1}^{(\theta)}\le\theta_0$ 
for~$\theta\le\theta_0$.
The bound \eqref{eins} implies that
\begin{equation}
\bigl|\E_{2n}^{(\theta)}(f(Q_{2n,1}^{(\theta)}))-\E_{2n}^{(\theta)}(f(Q_{n,1}^{(\theta_b)}))\bigr|
\le C'\Vert f\Vert_\infty e^{-\eta n},
\end{equation}
where~$C'<\infty$ and~$\eta>0$, while the bound \eqref{zwei} guarantees that~$\wE(f(Q_\infty))$ can be replaced by~$\wE(f(Q_{n,1}^{(\theta_b)}))$ with a similar error.
Then \eqref{Hbd} with~$2n$ replacing~$n$
boils down to the estimate of
\begin{equation}
\label{drei}
\Bigl|\E_{2n}^{(\theta)}\bigl(f(Q_{n,1}^{(\theta)})\bigr)-\wE\bigl(f(Q_{n,1}^{(\theta)})\bigr)\Bigr|.
\end{equation}
But, by Proposition~\ref{prop3.5}, the latter is bounded by~$A\Vert f\Vert_\infty e^{-\zeta n}$.
Combining all of the previous estimates, the claim follows.
\end{proofsect}

\begin{proofsect}{Proof of Corollary~\ref{cor3.7}}
We begin by showing that~$\zz(\rho)=\wBbbP(Q_\infty\ge1)$. Indeed,
we can use that~$Z_n(\theta)=0$ for~$\theta<1$ to compute
\begin{equation}
\begin{aligned}
\wE\bigl(Z_n(Q_\infty)\bigr) &=\wBbbP(Q_\infty\ge1)\, \E\otimes\wE
\Bigl(Z_{n-1}\bigl(X+{\textstyle\frac1b}Q_\infty\bigr)\,\Big|\,
Q_\infty\ge1\Bigr)
\\
&=\wBbbP(Q_\infty\ge1)\,\wE\bigl(Z_{n-1}(Q_\infty)\bigr)
=\dots=\wBbbP(Q_\infty\ge1)^{n+1},
\end{aligned}
\end{equation}
where we used Proposition~\ref{prop3.6} to derive the second
equality. From here~$\zz(\rho)=\wBbbP(Q_\infty\ge1)$ follows by
noting that~$\wBbbP(Q_\infty\ge1)Z_n(1)\le \wE(Z_n(Q_\infty))\le
Z_n(\theta_b)$ and applying Theorem~\ref{thm2.4}(1).

In order to prove the existence of the limit \eqref{psirho}, we first notice that
\begin{equation}
\frac{Z_{n+1}(\theta)}{Z_n(\theta)}
=\BbbP_{n+1}^{(\theta)}\bigl(Q_{n+1,1}^{(\theta)}\ge1\bigr).
\end{equation}
Next we claim that~$\BbbP_{n+1}^{(\theta)}(Q_{n+1,1}^{(\theta)}\ge1)-\zz(\rho)$, for~$\theta\ge1$, decays exponentially with~$n$. Indeed,
let~$\theta_0>\theta_b$ and~$\theta\in[1,\theta_0]$, pick~$k=\lfloor \frac n2\rfloor$, use~$Q_{k,1}^{(1)}\le Q_{n+1,1}^{(\theta)}\le
Q_{k,1}^{(\theta_b)}$ and apply Proposition~\ref{prop3.5},
to get
\begin{equation}
\wBbbP\bigl(Q_{k,1}^{(1)}\ge1\bigr)-\bar Ae^{-\zeta k} \le
\BbbP_{n+1}^{(\theta)}\bigl(Q_{n+1,1}^{(\theta)}\ge1\bigr)\le
\wBbbP\bigl(Q_{k,1}^{(\theta_0)} \ge1\bigr)+\bar Ae^{-\zeta k},
\end{equation}
where~$\bar A<\infty$ is proportional to~$A(\rho,\theta_0)$ from~\eqref{cbd}.
On the other hand, we clearly have
\begin{equation}
\wBbbP\bigl(Q_{k,1}^{(1)}\ge1\bigr)\le \wBbbP(Q_\infty\ge1)\le
\wBbbP\bigl(Q_{k,1}^{(\theta_0)}\ge1\bigr).
\end{equation}
But the right and left-hand sides of this inequality differ only
by~$\wBbbP(Q_{k,1}^{(1)}<1,\, Q_{k,1}^{(\theta_0)}\ge1)$,
which can be estimated as in \eqref{QQdiff} by a number tending to
zero exponentially fast as~$k\to\infty$. From here we have
\begin{equation}
\label{3.16}
\biggl|\frac{Z_{n+1}(\theta)}{Z_n(\theta)\zz(\rho)}-1\biggr|\le
A' e^{-\zeta' n},\qquad\theta\in[1,\theta_0],
\end{equation}
where~$A'=A'(\rho,\theta_0)<\infty$ and~$\zeta'=\zeta'(\rho)>0$. 
The uniformity of these estimates
is a consequence of the uniformity of the
bounds~$A<\infty$ and~$\zeta>0$ and that 
as in \eqref{QQdiff}.

The existence of the limit \eqref{psirho} for~$\theta\in[1,\theta_0]$ is a direct consequence of \eqref{3.16} and the~identity
\begin{equation}
\label{prdct}
\psi_\rho(\theta)=\lim_{n\to\infty}Z_n(\theta)\zz(\rho)^{-n}
=\lim_{n\to\infty}\prod_{k=0}^{n-1}\frac{Z_{k+1}(\theta)}{\zz(\rho)Z_k(\theta)}
=\prod_{k=0}^\infty\frac{Z_{k+1}(\theta)}{\zz(\rho)Z_k(\theta)},
\end{equation}
and the fact that the corresponding infinite product converges.
For~$\theta<1$ we have~$Z_n(\theta)=0$ and the limit exists trivially. To prove that~$\theta\mapsto\psi_\rho(\theta)$ is Lipschitz continuous for~$\theta\ge1$, 
we first note
that, by \eqref{three} and the result of Lemma~\ref{lemma5.8a},
\begin{equation}
\bigl|Z_n(\theta)-Z_n(\theta')\bigr|\le C |\theta-\theta'|\psi_\rho(\theta_0)
\zz(\rho)^{-n},\qquad
\theta,\theta'\in[1,\theta_0],
\end{equation}
where~$C=C(\rho,\theta_0)<\infty$ is
on sets~$\NN\subset\MMsharp$ with finitely many extreme points. From here
the bound in part (2) directly follows.

Let~$Z_n^{(\rho)}(\theta)$ denote explicitly that~$Z_n(\theta)$ is computed
using the underlying measure~$\rho$.
The continuity of~$\alpha\mapsto\psi_{\rho_\alpha}(\theta)$ then follows using three facts: First,~$\alpha\mapsto Z_n^{(\rho_\alpha)}(\theta)$,
being an expectation with respect to~$\rho_\alpha^n$,
is continuous. Second, by Theorem~\ref{thm2.4}(2),~$\alpha\mapsto\zz(\rho_\alpha)$ is also continuous. Third, the infinite product \eqref{prdct} converges uniformly in~$\alpha$.
\end{proofsect}

\vbox{
\section*{Acknowledgements}
\noindent
We wish to thank T.~Liggett and R.~Lyons for useful discussions.
The research of L.C.~was supported by the NSF under the grant DMS-9971016 and
by the NSA under the grant NSA-MDA~904-00-1-0050.
}

\end{document}